\documentclass[12pt,oneside]{amsart}
\newtheorem{thm}{Theorem}[section] 
\newtheorem{cor}[thm]{Corollary}
\newtheorem{lem}[thm]{Lemma}
\newtheorem{prop}[thm]{Proposition}
\theoremstyle{definition}
\theoremstyle{remark}
\newtheorem{rem}[thm]{Remark}
\theoremstyle{proof}

\numberwithin{equation}{section}
\newcommand{\norm}[1]{\left\Vert#1\right\Vert}
\newcommand{\snorm}[1]{\Vert#1\Vert}
\newcommand{\abs}[1]{\left\vert#1\right\vert}
\newcommand{\set}[1]{\left\{#1\right\}}
\newcommand{\brac}[1]{\left(#1\right)}
\newcommand{\scalar}[1]{\left \langle #1 \right \rangle}
\newcommand{\abscalar}[1]{\abs{\scalar{#1}}}
\newcommand{\Real}{\mathbb{R}}

\newcommand{\To}{\longrightarrow}

\newcommand{\E}{\mathcal{E}}
\newcommand{\BP}{BP}
\renewcommand{\MR}{\widetilde{M}} 
\newcommand{\M}{\mathcal{M}}
\newcommand{\I}{\mathcal{I}}
\renewcommand{\L}{\mathcal{L}}
\newcommand{\Lk}{\mathcal{L}_k}

\newcommand{\Vol}[1]{\textnormal{Vol} (#1)}
\newcommand{\sVolRatio}[2]{\brac{\Vol{#1}/\Vol{#2}}^{1/n}}
\newcommand{\VolRatio}[2]{\brac{\frac{\Vol{#1}}{\Vol{#2}}}^{1/n}}
\newcommand{\definition}[1]{\medskip\noindent\textbf{Definition.}#1\medskip}
\begin{document}


%

\title{Dual Mixed Volumes and the Slicing Problem}
\author{Emanuel Milman\\ \bigskip \fontfamily{cmr} \fontseries{m} \fontshape{sc} \fontsize{10}{0} \selectfont
Department of Mathematics, \\
The Weizmann Institute of Science, \\Rehovot 76100, Israel. \\
\medskip E-Mail: emanuel\_milman@hotmail.com.}
\thanks{Supported in part by BSF and ISF}
\begin{abstract}
We develop a technique using dual mixed-volumes to study the
isotropic constants of some classes of spaces. In particular, we
recover, strengthen and generalize results of Ball and Junge
concerning the isotropic constants of subspaces and quotients of
$L_p$ and related spaces. An extension of these results to
negative values of $p$ is also obtained, using generalized
intersection-bodies. In particular, we show that the isotropic
constant of a convex body which is contained in an
intersection-body is bounded (up to a constant) by the ratio
between the latter's mean-radius and the former's volume-radius.
We also show how type or cotype 2 may be used to easily prove
inequalities on any isotropic measure.

\end{abstract}



\maketitle


\section{Introduction}

The main purpose of this note is to provide new types of bounds on
a convex body's isotropic constant, by means of dual mixed-volumes
with different families of bodies.

A centrally symmetric convex body $K$ in $\Real^n$ is said to be
in isotropic position if $\int_K \scalar{x,\theta}^2 dx$ is
constant for all $\theta \in S^{n-1}$, the Euclidean unit sphere.
If in addition $K$ is of volume 1, then its isotropic constant is
defined to be the $L_K$ satisfying $\int_K \scalar{x,\theta}^2 dx
= L_K^2$ for all $\theta \in S^{n-1}$.
It is easy to see that every body may be brought to isotropic
position using an affine transformation, and that the isotropic
position is unique modulo rotations and homothety
(\cite{Milman-Pajor-LK}). Hence, for a general centrally
symmetric convex body $K$ we shall denote by $L_K$ the isotropic
constant of $K$ in its isotropic position of volume 1.

A famous problem, commonly known as the Slicing Problem, asks
whether $L_K$ is bounded from above by a universal constant
independent of $n$, for all centrally symmetric convex bodies $K$
in $\Real^n$. This was first posed in an equivalent form by J.
Bourgain, who asked whether every centrally symmetric convex body
of volume 1, has an $n-1$ dimensional section whose volume is
bounded from below by some universal constant. This is known to
be true for several families of bodies, such as sections of $L_1$,
projection bodies and 1-unconditional bodies (see
\cite{Milman-Pajor-LK},\cite{Ball-weak-GL} or below). The best
general bound is due to Bourgain, who showed in \cite{Bourgain-LK}
that $L_K \leq C n^{1/4} \log(1+n)$. Recently, the general
problem has been reduced to the case that $K$ has finite
volume-ratio (\cite{BKM-symmetrizations}).

The main idea of this note is to compare a general convex body
$K$ (or its polar) with a less general body $L$ chosen from a
specific family, and thus gain some knowledge on its isotropic
constant. We shall consider two main families: unit-balls of
$n$-dimensional subspaces of $L_p$, denoted $SL_p^n$, and
$k$-Busemann-Petty bodies, denoted $\BP_k^n$, which are a
generalization of intersection bodies (the class $\BP_1^n$)
introduced by Zhang in \cite{Zhang-Gen-BP} (there they are
referred to as "generalized $(n-k)$-intersection bodies", see
Section \ref{sec:definitions} for definitions). The body $L$ may
not be necessarily convex, but we will assume that it is a
centrally symmetric star-body, defined by a continuous radial
function $\rho_L(\theta) = \max \set{r \geq 0 | r\theta \in L }$
for $\theta \in S^{n-1}$. Our main tool for comparing two
star-bodies will be the dual mixed-volume of order $p$, defined
in Section \ref{sec:definitions}, which was first introduced by
Lutwak in \cite{Lutwak-dual-mixed-volumes}.

We will require a few more notations. Let $\abs{x}$ denote the
standard Euclidean norm of $x \in \Real^n$, let $D_n$ denote the
Euclidean unit ball and let $\sigma$ denote the Haar probability
measure on $S^{n-1}$. Let $SL(n)$ denote the group of volume
preserving linear transformations in $\Real^n$, and let $\Vol{B}$
denote the Lebesgue measure of the set $B \subset \Real^n$ in its
affine hull. Let $K^\circ$ denote the polar body to a convex body
$K$.

An equivalent characterization of the isotropic position
(\cite{Milman-Pajor-LK}) states that it is the position which
minimizes the expression $\int_K \abs{x}^2 dx$, in which case the
latter is equal to $n L_K^2$ if $\Vol{K}=1$. By comparing with
the value of this expression in a position for which the
circumradius $a(K)$ of $K$ is minimal, we immediately get the
bound $L_K \leq a(K) / \sqrt{n}$. Equivalently, making this
invariant to change of position or normalization, we get the
following well known elementary bound on $L_K$ in terms of the
outer volume-ratio of $K$:
\[
L_K \leq C \inf \set{\left . \VolRatio{\E}{K} \right | \; K
\subset \E \; , \; \E \in SL_2^n},
\]
where $SL_2^n$ is just the class of all ellipsoids in $\Real^n$.
This was generalized in \cite{Ball-weak-GL} by K. Ball as follows:

\smallskip
\noindent\textbf{Theorem (Ball).}
\begin{equation} \label{eq:Ball}
L_K \leq C \inf \set{ \left . \VolRatio{L}{K} \right | \: K
\subset L \; , \; L \in SL_1^n}.
\end{equation}
In fact, Ball showed that the expression on the right is
equivalent (up to universal constants) to the so-called weak
right-hand Gordon-Lewis constant $wrgl_2(X^*_K)$ of the Banach
space $X^*_K$ whose unit ball is the polar of $K$. Ball showed
that $wrgl_2(X^*)$ is majorized (up to a constant) by $gl_2(X)$,
the Gordon-Lewis constant of $X$, and hence $L_K$ is bounded for
spaces $X_K$ with uniformly bounded $gl_2$ constants. These
include subspaces of $L_p$ for $1\leq p \leq 2$, quotients of
$L_q$ for $2\leq q \leq \infty$, and spaces with a
1-unconditional basis (the latter were first shown to have a
bounded isotropic constant by Bourgain). A complementary result
was obtained in \cite{Junge-slicing-problem-for-quotients-of-Lp}
by Junge, who showed the following (this is not explicit in his
formulation but follows from the proof):

\smallskip
\noindent\textbf{Theorem (Junge). }
\begin{equation} \label{eq:Junge}
L_K \leq C \inf \set{ \left . \sqrt{p} \; q \VolRatio{L}{K} \:
\right
| \; \begin{array}{c} K \subset L \; , \; L \in SQL_p^n \; , \\
1<p<\infty \; , \; 1/p + 1/q = 1 \end{array} },
\end{equation}
where $SQL_p^n$ is the class of all unit-balls of $n$-dimensional
subspaces of quotients of $L_p$, and $q = p^*$ is the conjugate
exponent to $p$. In fact, Junge showed that $L_p$ may be replaced
by any Banach space $X$ with bounded $gl_2(X)$ such that $X$ has
finite type, in which case $\sqrt{p} \; q$ above should be
replaced by some constant depending on $X$.

\medskip
As evident from their more general formulations, the results of
Ball and Junge described above make heavy use of non-trivial
Functional Analysis and Operator Theory, and as a result the
geometric intuition behind the Slicing Problem is substantially
lost. Of course, this is to be expected if the conditions on the
space $X_K$ are formulated using Operator Theory notions, such as
(variants of) the Gordon-Lewis property. But for classical spaces
such as subspaces or quotients of $L_p$, one may hope to simplify
the approach, derive better bounds on $L_K$, and unify Ball and
Junge's results into a single framework. Using an elementary
argument, geometric in nature, we show the following
generalizations of (\ref{eq:Ball}) and partial strengthening of
(\ref{eq:Junge}) (the term ``partial" refers to the fact that we
restrict $L$ to the class $SL_p^n$ or $QL_q^n$ defined below),
for a convex isotropic body $K$ with $\Vol{K} = \Vol{D_n}$:

\smallskip
\noindent\textbf{Theorem 1.}
\[
L_K \leq C \; \inf \set{ \left . \frac{\sqrt{p_0}}{M_p(L)} \right
| \: K \subset L \; , \; L \in SL_p^n \; , \; p\geq 0},
\]
where $p_0 = \max(1,\min(p,n))$, $M_p(L) = \brac{\int_{S^{n-1}}
\norm{x}^p_L d\sigma(x)}^{1/p}$ for \mbox{$p>0$}, and by passing
to the limit, $M_0(L) = \exp \brac{\int_{S^{n-1}} \log \norm{x}_L
d\sigma(x)}$.

\smallskip
\noindent\textbf{Theorem 1'.}
\[
L_K \leq C \frac{T_2(X_K)}{M_2(K)},
\]
where $T_2(X_K)$ is the (Gaussian) type-2 constant of $X_K$.

\smallskip
\noindent\textbf{Theorem 2.}
\[
L_K \leq C \; \inf \set{ \left . \Lk \MR_k(L) \right | \: K
\subset L \; , \; L \in \BP_k^n \; , \; k = 1,\ldots,n-1 },
\]
where $\mathcal{L}_k$ denotes the maximal isotropic constant of
centrally symmetric convex bodies in $\Real^k$ and $\MR_k(L) =
\brac{\int_{S^{n-1}} \rho_L(x)^k d\sigma(x)}^{1/k}$. We emphasize
again that $\BP_1^n$ is exactly the class of intersection bodies.

\medskip \noindent Indeed, these are all generalizations of
(\ref{eq:Ball}) and (\ref{eq:Junge}), since by passing to polar
coordinates and applying Jensen's inequality (for $p,k>0$):
\begin{equation} \label{eq:Jensen}
\frac{1}{M_p(L)} \leq \MR_k(L) \leq \VolRatio{L}{D_n},
\end{equation}
and since $T_2(X_K) \leq C \sqrt{p}$ by Kahane's inequality when
$K \in SL_p^n$ for $p \geq 2$. This also applies to Theorem 2,
since any $K \in SL^n_p$ for $0<p\leq 2$ (and in particular
$p=1$) is an intersection body (see
\cite{Koldobsky-intersection-bodies}), and hence a
$k$-Busemann-Petty body for all $k\geq 1$
(\cite{Grinberg-Zhang},\cite{EMilman-Generalized-Intersection-Bodies}).

\medskip

We also have the following dual counterparts to Theorems 1 and 2,
for a convex isotropic body $K$ with $\Vol{K} = \Vol{D_n}$:

\smallskip
\noindent\textbf{Theorem 3.}
\smallskip
\[
L_K \leq C \; \inf \set{ \sqrt{p_0} \; M^*_p(T(L)) \; \left | \:
\begin{array}{c} K \subset L \; , \; L \in QL_q^n \; , \; T \in
SL(n) \; , \\  1\leq q\leq\infty \; , \; 1/p + 1/q = 1
\end{array} \right . },
\]

\noindent where $QL_q^n$ is the class of all unit-balls of
$n$-dimensional quotients of $L_q$, $p_0$ is defined as above for
$p = q^*$, and $M^*_p(G) = M_p(G^\circ)$.

\medskip
\noindent This is indeed a (partial) strengthening of
(\ref{eq:Junge}), since by Lemma \ref{lem:M_p-bound} (see also
the Mean Norm Corollary below), there exists a position $T \in
SL(n)$ of $L\in QL_q^n$ such that:
\[
M^*_p(T(L)) \leq C \sqrt{p_0} \VolRatio{L}{D_n}.
\]
It is also interesting to note that the proof of Theorem 3,
although derived independently, closely resembles Bourgain's
proof that $L_K \leq C n^{1/4} \log(1+n)$.

\medskip
\smallskip
\noindent\textbf{Theorem 4.}
\smallskip
\[
L_K \leq C \inf \set{\left . \frac{\L_{2k}^2}{\MR_k(T(L))}  \right
|
\begin{array}{c} L \subset K^\circ \; , \; L \in \BP_k^n \; , \\ T \in SL(n) \; , \; k =
1,\ldots,\lfloor n/3 \rfloor \end{array} }.
\]

\noindent Using an analogue of Lemma \ref{lem:M_p-bound} (stated
in the Mean Radius Corollary below), we may deduce the following
bound on $L_K$ for polars of bodies in $C\BP_k^n$, the class of
\emph{convex} $k$-Busemann-Petty bodies:

\[
L_K \leq C \inf \set{\left . \L_{2k}^2 \Lk \VolRatio{L^\circ}{K}
\right | \begin{array}{c} K \subset L^\circ \; , \; L \in
C\BP_k^n \; , \\ k = 1,\ldots,\lfloor n/3 \rfloor \end{array} }.
\]

\medskip

Since Jensen's inequality (\ref{eq:Jensen}) is usually strict, it
is not hard to construct examples for which Theorem 1
asymptotically out-performs Junge's bound. Indeed, for $K =
[-1,1]^n$, it is well known (see Section \ref{sec:applications})
that $K$ is isomorphic to a body $L \in SL_p^n$, for $p = \log
n$. Junge's bound therefore implies $L_K \leq C \sqrt{\log n}$,
while Theorem 1 gives $L_K \leq C$, since $M_p(L) \simeq M_p(K)
\simeq \sqrt{\log n} \sVolRatio{K}{D_n}$.

\smallskip

As mentioned above, Theorems 1 and 3 imply, in particular, that
$L_K \leq C \sqrt{p}\;$ for $K \in SL_p^n$ and $p \geq 1$, and
$L_K \leq q^*$ for $K \in QL_q^n$ and $q > 1$. We note that this
is not contained in Junge's result (\ref{eq:Junge}). The strength
of (\ref{eq:Junge}) is that it applies simultaneously to all
subspaces of quotients of $L_p$, which our method does not
handle. Ironically, this is also its drawback, if one is
interested in proper subspaces or quotients only: it gives the
same bound on $L_K$ in either case. Therefore, one cannot hope to
have a good bound for $SL_p^n$ with $1 \leq p < 2$ ($QL_q^n$ with
$q>2$) without solving the Slicing Problem, because this would
imply the same bound for $QL_p^n$ ($SL_q^n$) in that range, which
already contain all convex bodies. To fill the bound for $SL_p^n$
with $1 \leq p < 2$ ($QL_q^n$ with $q>2$), one needs to use Ball's
result in its general form (or simply use (\ref{eq:Ball})
combined with the fact that $SL_p^n \subset SL_1^n$ for $1 \leq p
\leq 2$; by duality $QL_q^n \subset QL_\infty^n$ for $q \geq 2$,
implying that the bodies in $QL_q^n$ have finite outer
volume-ratio as projection bodies). We therefore see that Theorems
1 and 3 combine the ranges $1 \leq p < 2$ and $p\geq 2$ into a
single framework.

\smallskip

Evidently, Theorem 1' has a somewhat different flavor, and indeed
its proof is totally different from the proofs of the other
Theorems. The proof is based on a simple yet effective framework
for combining isotropic measures with type and cotype 2, which is
introduced in Section \ref{sec:type-2} (this section may be read
independently from the rest of this note). This framework also
enables us to easily recover several known lemmas on John's
maximal volume ellipsoid position (originally proved using
Operator Theory techniques), which we use in the proof of Lemma
\ref{lem:M_p-bound} (mentioned above). We remark that Theorem 1'
also follows from the work in \cite{BMMP} but in a more
complicated manner.

\smallskip

The other Theorems are all proved using another technique,
involving dual mixed-volumes.
Theorems 1 and 3 are proved in Section \ref{sec:SectionsOfLp}, and
Theorems 2 and 4 are proved in Section \ref{sec:k-BP-bodies}. In
Section \ref{sec:applications}, we give several corollaries of
our main Theorems, some of which are mentioned below.

\smallskip

Using the known fact that $L_K$ is always bounded from below,
Theorems 1',1 and 2, immediately yield the following useful
corollary, for an isotropic convex body $K$ with
$\Vol{K}=\Vol{D_n}$:

\smallskip
\noindent\textbf{Mean Norm/Radius Corollary.}\emph{
\begin{enumerate}
\item
$M_2(K) \leq C T_2(X_K)$. \vspace{3pt}
\item
If $K \in SL_p^n$ ($p>0$), then $M_p(K) \leq C \sqrt{p_0}$.
\item
If $K \in \BP_k^n$ ($k=1,\ldots,n-1$), then $\MR_k(K) \geq C /
\Lk$.
\end{enumerate}
}

\noindent Jensen's inequality in (\ref{eq:Jensen}) shows that
these bounds are tight (to within a constant) for $p,k,T_2(X_K)
\leq C$. One should also keep in mind that if $K^\circ$ is in
isotropic position, this corollary is applicable to $K^\circ$,
providing different inequalities.

\medskip

In addition, although this is a direct consequence of the
extended formulation of Junge's Theorem (and also of Theorems 1
and 3), the following corollary about a centrally symmetric convex
polytope $P$ is worth explicit stating:
\smallskip

\noindent\textbf{Polytope Corollary.} \emph{
\begin{enumerate}
\item
If $P$ has $2m$ facets then $L_P \leq C \sqrt{\log(1+m)}$.
\item
If $P$ has $2m$ vertices then $L_P \leq C \log(1+m)$.
\end{enumerate}
} \noindent In particular, this implies that Gluskin's
probabilistic construction in \cite{Gluskin-Diameter} of two
convex bodies $K_1$ and $K_2$ with Banach-Mazur distance of order
$n$, satisfies $L_{K_1},L_{K_2} \leq C \log(1+n)$.

\medskip

Theorem 2 should be understood as a partial complimentary result
to Theorem 1. The reason for this may be better explained, if we
first consider a second generalization of intersection bodies,
introduced by Koldobsky in \cite{Koldobsky-I-equal-BP}. We shall
call these bodies \emph{$k$-intersection bodies} and denote this
class of bodies by $\I_k^n$. It was shown in
\cite{Koldobsky-I-equal-BP} that $\BP_k^n \subset \I_k^n$, and
the question of whether $\BP_k^n = \I_k^n$ remains open (see
\cite{EMilman-Generalized-Intersection-Bodies} for an account of
recent progress in this direction). The class $\I_k^n$ satisfies
a certain characterization of being embedded in $L_p$, which has
been continued analytically to the negative value $p=-k$, so in
some sense $\I_k^n = SL_{-k}^n$. Therefore, in some sense,
$\BP_k^n \subset SL_{-k}^n$, hence our initial remark.

The class of star-bodies $\BP_k^n$ seems at first glance a
non-natural object to work with when studying convex bodies.
Nevertheless, we describe in Section \ref{sec:applications}
several potential ways in which this object may be harnessed to
our advantage.

%

\medskip

\noindent \textbf{Acknowledgments.} I would like to deeply thank
my supervisor Prof. Gideon Schechtman for many informative
discussions, and especially for believing in me and allowing me to
pursue my interests. I would also like to thank the referee for
many helpful remarks.

\section{Definitions and Notations} \label{sec:definitions}

A convex body $K$ will always refer to a compact, convex set in
$\Real^n$ with non-empty interior. We will always assume that the
bodies in question are centrally symmetric, i.e. $K=-K$. The
equivalence between convex bodies and norms in $\Real^n$ is well
known, with the correspondence $\norm{x}_K = \min \set{t>0 | x/t
\in K}$. The associated normed space $(\Real^n,\norm{\cdot}_K)$
will be denoted by $X_K$. The dual norm is defined as
$\norm{x}^*_K = \sup_{y \in K} \abs{\scalar{x,y}}$, and its
associated unit-ball is called the polar body to $K$, and denoted
$K^\circ$. The dual normed space $(\Real^n,\norm{\cdot}^*_K)$ is
denoted by $X^*_K$ ($= X_{K^\circ}$). We will say that a
convex-body $K$ is 1-unconditional, or simply unconditional, with
respect to the given Euclidean structure (which we always assume
to be fixed), if $(x_1,\ldots,x_n) \in K$ implies $(\pm x_1,
\ldots, \pm x_n) \in K$ for all possible sign assignments.

We will also work with general star-bodies $L$, which are
star-shaped bodies, meaning that $tL \subset L$ for all $t \in
[0,1]$, with the additional requirement that their radial
function $\rho_L$ is a continuous function on $S^{n-1}$. The
radius of $L$ in direction $\theta \in S^{n-1}$ is defined as
$\rho_L(\theta) = \max \{\mbox{$r \geq 0$} \; | \; r\theta \in L
\}$. For a general star-body $L$, we define its Minkowski
functional $\norm{x}_L$ in the same manner as for a convex body
(so $\norm{x}_L$ is no longer necessarily a norm). Obviously,
$\rho_L(\theta) = 1/\norm{\theta}_L$ for all $\theta \in S^{n-1}$.

By identifying between a star-body and its radial function, a
natural metric arises on the space of star-bodies. The radial
metric, denoted by $d_r$, is defined as:
\[
d_r(L_1,L_2) = \sup_{\theta \in S^{n-1}} \abs{\rho_{L_1}(\theta)-
\rho_{L_2}(\theta)}.
\]

As mentioned in the Introduction, our main tool for comparing two
star-bodies $L_1$ and $L_2$ will be the dual mixed-volume of
order $p \in \Real$, introduced by Lutwak in
\cite{Lutwak-dual-mixed-volumes} (see also
\cite{Lutwak-intersection-bodies}), and defined as:
\[
\widetilde{V}_p(L_1,L_2) = \frac{1}{n} \int_{S^{n-1}}
\rho_{L_1}(x)^p \rho_{L_2}(x)^{n-p} \, dx
\]
(note that the integration is w.r.t. the Lebesgue measure on
$S^{n-1}$). By polar integration, it is obvious that
$\widetilde{V}_p(L,L) = \Vol{L}$ for all $p$. We will also use
the following useful property of dual mixed-volumes (see
\cite{Lutwak-intersection-bodies}):
\begin{equation} \label{eq:dmv-property}
\widetilde{V}_p(T(L_1),T(L_2)) = \widetilde{V}_p(L_1,L_2),
\end{equation}
for any $T \in SL(n)$ and $p \in \Real$. We also constantly use
the well known formula for the volume of the Euclidean unit ball
$D_n$:
\begin{equation} \label{eq:ball-volume-formula}
\Vol{D_n} = \frac{\pi^{n/2}}{\Gamma(n/2+1)}.
\end{equation}

\smallskip
Several useful notations for a star-body $L$ will be used. For
$p>0$, the $p$-th mean-norm, denoted by $M_p(L)$, is defined as:
\[
M_p(L) = \brac{\int_{ S^{n-1} } \norm{x}^p_L d\sigma(x)}^{1/p}.
\]
Passing to the limit as $p \rightarrow 0$, we define $M_0(L) =
\exp\brac{\int_{ S^{n-1} } \log \norm{x}_L d\sigma(x)}$. We will
define the mean-norm as $M(L) = M_1(L)$. The $p$-th mean-width,
denoted $M^*_p(L)$, is defined as $M^*_p(L) = M_p(L^\circ)$, and
as usual, the mean-width is defined as $M^*(L) = M^*_1(L)$. The
$p$-th mean-radius, denoted by $\MR_p(L)$, is defined as:
\[
\MR_p(L) = \brac{\int_{ S^{n-1} } \rho_L(x)^p d\sigma(x)}^{1/p}.
\]
We will define the mean-radius as $\MR(L) = \MR_1(L)$. The
minimal $a,b>0$ for which $1/a \abs{x} \leq \norm{x}_L \leq b
\abs{x}$, will be denoted by $a(L)$ and $b(L)$, respectively.
Geometrically, $a(L)$ and $1/b(L)$ are the radii of the
circumscribing and inscribed Euclidean balls of $L$, respectively.
the The expression $\sVolRatio{L}{D_n}$ will be referred to as the
volume-radius of $L$. The infimum of $\sVolRatio{L}{\E}$ over all
ellipsoids $\E$ contained in $L$ is called the volume-ratio of
$L$. Similarly, the infimum of $\sVolRatio{\E}{L}$ over all
ellipsoids $\E$ containing $L$ is called the outer volume-ratio
of $L$. A position of a body $L$ is a volume preserving linear
image of $L$, i.e. $T(L)$ for $T \in SL(n)$.

Going back to convex bodies and normed spaces, we now define the
(Gaussian) type and cotype 2 constants of a normed space
$X=(\Real^n,\norm{\cdot})$. The (Gaussian) type-2 constant of $X$,
denoted $T_2(X)$, is the minimal $T>0$ for which:
\[
\brac{\int_\Omega \snorm{\sum_{i=1}^m g_i(\omega) x_i}^2
d\omega}^{1/2} \leq T \brac{\sum_{i=1}^m \norm{x_i}^2}^{1/2}
\]
for any $m \geq 1$ and any $x_1,\ldots,x_m \in X$, where
$g_1,\ldots,g_m$ are independent real-valued standard Gaussian
r.v.'s on a common probability space $(\Omega,d\omega)$.
Similarly, the (Gaussian) cotype-2 constant of $X$, denoted
$C_2(X)$, is the minimal $C>0$ for which:
\[
\brac{\int_\Omega \snorm{\sum_{i=1}^m g_i(\omega) x_i}^2
d\omega}^{1/2} \geq 1/C \brac{\sum_{i=1}^m \norm{x_i}^2}^{1/2}
\]
for any $m \geq 1$ and $x_1,\ldots,x_m \in X$. We will not
distinguish between the Gaussian and the Rademacher type (cotype)
2 constants, since it is well known that the former constant is
always majorated by the latter one (e.g.
\cite{Milman-Schechtman-Book}), and all our results will involve
upper bounds in terms of the Gaussian type (cotype) 2.

We will often identify between a normed space and its unit-ball.
In particular, for the infinite dimensional Banach space $L_p =
L_p([0,1],dx)$, whenever the expression "sections of $L_p$" is
used, we will mean sections of its unit-ball. And when the
expression "quotients of $L_p$" is used, we might refer to the
unit-balls of these quotient spaces.

Throughout the paper, all constants used will be universal,
independent of all other parameters, and in particular,
independent of $n$. We reserve $C,C',C_1,C_2$ to denote these
constants, which may take different values on separate instances.
We will write $A \simeq B$ to signify that $C_1 A \leq B \leq C_2
A$ with universal constants $C_1,C_2>0$.

\medskip

For the results of Sections \ref{sec:k-BP-bodies} and
\ref{sec:applications}, we shall need to define the class of
$k$-Busemann-Petty bodies, introduced by Zhang in
\cite{Zhang-Gen-BP} (there they are referred to as "generalized
$(n-k)$-intersection bodies"). These bodies represent a
generalization of the notion of an intersection body. For
completeness, we give the appropriate definitions below.

\definition{
A star body $K$ is said to be an
\emph{intersection body of a star body} $L$, if $\rho_K(\theta) =
\Vol{L \cap \theta^\perp}$ for every $\theta \in S^{n-1}$.
$K$ is said to be an \emph{intersection body}, if it is the limit
in the radial metric $d_r$ of intersection bodies $\{K_i\}$ of
star bodies $\{L_i\}$. This is equivalent (e.g.
\cite{Lutwak-intersection-bodies}, \cite{Gardner-BP-5dim}) to
$\rho_K = R^*(d\mu)$, where $\mu$ is a non-negative Borel measure
on $S^{n-1}$, $R^*$ is the dual transform (as in
(\ref{eq:duality111})) to the Spherical Radon Transform
$R:C(S^{n-1}) \rightarrow C(S^{n-1})$, which is defined for $f\in
C(S^{n-1})$ as:
\[
R(f)(\theta) = \int_{S^{n-1} \cap \theta^\perp} f(\xi)
d\sigma_{n-1}(\xi),
\]
where $\sigma_{n-1}$ the Haar probability measure on $S^{n-2}$
(and we have identified $S^{n-2}$ with $S^{n-1} \cap
\theta^\perp$).
}

Let $G(n,m)$ denote the Grassmann manifold of all $m$-dimensional
linear subspaces of $\Real^n$. Generalizing the Spherical Radon
Transform is the $m$-dimensional Spherical Radon Transform $R_m$,
acting on spaces of continuous functions as follows:
\begin{eqnarray}
\nonumber R_m: C(S^{n-1}) \To C(G(n,m)) \\
\nonumber R_m(f) (E) = \int_{S^{n-1}\cap E} f(\theta)
d\sigma_m(\theta) ,
\end{eqnarray}
where $\sigma_m$ is the Haar probability measure on $S^{m-1}$ (and
we have identified $S^{m-1}$ with $S^{n-1} \cap E$). Notice that
for a star-body $L$ in $\Real^n$:
\[
R_m(\rho^m_L)(E) =  \Vol{L \cap E} / \Vol{D_m} \;\; \forall E \in
G(n,m).
\]
The dual transform is defined on spaces of \emph{signed} Borel
measures $\M$ by:
\begin{eqnarray}
\label{eq:duality111} & R_m^*: \M(G(n,m)) \To \M(S^{n-1}) & \\
\nonumber & \int_{S^{n-1}} f R_m^*(d\mu) = \int_{G(n,m)} R_m(f)
d\mu & \forall f \in C(S^{n-1}),
\end{eqnarray}
and for a measure $\mu$ with continuous density $g$, the
transform may be explicitly written in terms of $g$ (see
\cite{Zhang-Gen-BP}):
\begin{eqnarray}
\nonumber R_m^* g (\theta) = \int_{\theta \in E \in G(n,m)} g(E)
d\nu_{m}(E) ,
\end{eqnarray}
where $\nu_{m}$ is the Haar probability measure on $G(n-1,m-1)$.

\definition{
A star body $K$ is said to be a \emph{$k$-Busemann-Petty body} if
$\rho_K^k = R_{n-k}^*(d\mu)$, where $\mu$ is a non-negative Borel
measure on $G(n,n-k)$. We shall denote the class of such bodies
by $\BP_k^n$. }

Choosing $k=1$, for which $G(n,n-1)$ is isometric to $S^{n-1} /
Z_2$ by mapping $H$ to $S^{n-1}\cap H^\perp$, and noticing that
$R$ is equivalent to $R_{n-1}$ under this map, we see that
$\BP_1^n$ is exactly the class of intersection bodies.

\medskip

To conclude this section, we mention that we always work with the
radial metric topology on the space of star-bodies. Equivalently,
we always work with the maximum norm on the space of continuous
functions on $S^{n-1}$. So whenever an expression of the
following form appears:
\[
f = \int f_\alpha d\mu(\alpha),
\]
where $f$ and $\set{f_\alpha}$ are continuous functions on
$S^{n-1}$, the convergence of the integral should be understood in
the maximum norm.



\section{Combining Isotropic Measures with Type / Cotype 2}
\label{sec:type-2}

In this section we introduce a very simple yet effective
framework, which demonstrates how to utilize isotropic measures
associated with a convex body $K$, to give bounds on $M_2(K)$ and
$M^*_2(K)$ in terms of the type-2 and cotype-2 constants of $X_K$
and $X^*_K$. As an immediate corollary, we revive a couple of
known (yet partially forgotten) lemmas on John's maximal volume
ellipsoid position, one of which will be used in Section
\ref{sec:SectionsOfLp} to improve the bound on the isotropic
constant of quotients of $L_q$. Another immediate corollary of
this framework is that $L_K$ is always bounded by $T_2(X_K)$.

\medskip

Recall that a Borel measure $\mu$ on $\Real^n$ is said to be
\emph{isotropic} if:
\[
\int_{\Real^n} \scalar{x,\theta}^2 d\mu(x) = \abs{\theta}^2
\;\;\; \forall \theta \in \Real^n.
\]
This is easily seen to be equivalent to:
\[
\int_{\Real^n} \scalar{x,\theta_1} \scalar{x,\theta_2} d\mu(x) =
\scalar{\theta_1,\theta_2} \;\;\; \forall \theta_1,\theta_2 \in
\Real^n.
\]
The main point of this section is the following easy yet useful
observation:
\begin{lem} \label{lem:isotropic_is_gaussian}
Let $v_i \in \Real^n$ and $\lambda_i > 0$, for $i=1,\ldots,m$, be
such that $\mu = \sum_{i=1}^m \lambda_i \delta_{v_i}$ is an
isotropic measure. Let $\set{g_i}_{i=1}^m$ be a sequence of
independent real-valued standard Gaussian r.v.'s, and define the
r.v. $\Lambda_\mu$ as:
\begin{equation} \label{eq:Lambda}
\Lambda_\mu = \sum_{i=1}^m g_i \sqrt{\lambda_i} v_i.
\end{equation}
Then $\Lambda_\mu$ is an $n$-dimensional standard Gaussian.
\end{lem}
\begin{proof}
Obviously $\Lambda_\mu$ is a zero mean Gaussian r.v., so it
remains to show that its correlation matrix is the identity.
Indeed, from the independence of the $g_i$'s and the isotropicity
of $\mu$:
\begin{eqnarray}
\nonumber E\brac{\scalar{\Lambda_\mu,\theta_1}
\scalar{\Lambda_\mu,\theta_2}} = E \brac{ \sum_{i,j=1}^m g_i g_j
\sqrt{\lambda_i} \sqrt{\lambda_j}  \scalar{v_i,\theta_1}
\scalar{v_j,\theta_2} } = \\
\nonumber E \brac{ \sum_{i=1}^m g_i^2 \lambda_i
\scalar{v_i,\theta_1} \scalar{v_i,\theta_2} } = \sum_{i=1}^m
\lambda_i \scalar{v_i,\theta_1} \scalar{v_i,\theta_2} =
\scalar{\theta_1,\theta_2}.
\end{eqnarray}
\end{proof}

By taking the Fourier transform of the densities on both sides of
(\ref{eq:Lambda}), or by projecting them onto an arbitrary
direction, we get:
\[
\exp(-\abs{x}^2) = \Pi_{i=1}^m
\brac{\exp(-\scalar{x,v_i}^2)}^{\lambda_i}.
\]
This formulation, which is easy to check directly, has been used
by many authors 
(e.g. \cite{Schechtman-Schlumprecht}, \cite{BallUseBrascampLieb}),
mostly with connection to John's decomposition of the identity.
The advantage of Lemma \ref{lem:isotropic_is_gaussian} is that we
may work directly on the Gaussian r.v.'s and use type and cotype
estimates on $\norm{\Lambda_\mu}$, as summarized in the following
Proposition.


\begin{prop} \label{prop:isotropic-prop}
Let $K$ denote a convex body and let $\mu$ be any finite,
compactly supported, isotropic measure. Then:
\[
\frac{1}{C_2(X_K)} \brac{\int \norm{x}^2_K d\mu(x)}^{1/2} \leq
\sqrt{n} M_2(K) \leq T_2(X_K) \brac{\int \norm{x}^2_K
d\mu(x)}^{1/2}
\]
\end{prop}
\begin{proof}
First, assume that $\mu$ is a discrete isotropic measure
supported on finitely many points, of the form $\mu = \sum_{i=1}^m
\lambda_i \delta_{v_i}$.
Then by Lemma \ref{lem:isotropic_is_gaussian}, denoting
$\set{g_i}_{i=1}^m$ and $\set{g'_i}_{i=1}^n$ two sequences of
independent standard Gaussian r.v.'s on a common probability
space $(\Omega,d\omega)$, we have:
\begin{eqnarray}
\nonumber \int_\Omega \snorm{\sum_{i=1}^m g_i(\omega)
\sqrt{\lambda_i} v_i}^2_K d\omega & = & \int_{\Omega}
\snorm{\sum_{i=1}^n g'_i(\omega) e_i}^2_K d\omega = \\
\nonumber \frac{1}{(2\pi)^{n/2}}\int_{\Real^n} \norm{x}^2_K
e^{-\abs{x}^2/2} dx & = & \frac{\int_0^\infty e^{-r^2/2}r^{n+1}
dr}{(2\pi)^{n/2}} \int_{S^{n-1}} \norm{\theta}^2_K d\theta = n
M_2(K)^2,
\end{eqnarray}
where the last equality is a standard calculation (e.g.
\cite{Milman-Schechtman-Book}). But on the other hand, using the
type-2 condition on $X_K$, we see that the initial expression on
the left is bounded from above by:
\[
T_2(X_K)^2 \sum_{i=1}^m \snorm{ \sqrt{\lambda_i} v_i}^2_K  =
T_2(X_K)^2 \int\norm{x}^2_K d\mu(x).
\]
Taking square root, the type-2 upper bound follows for a discrete
measure $\mu$, and the cotype-2 lower bound follows similarly.

When $\mu$ is a general isotropic measure, we approximate $\mu$
by a series of discrete (not necessarily isotropic) measures
$\mu_\epsilon = \sum_{i=1}^{m_\epsilon} \lambda^\epsilon_i
\delta_{v^\epsilon_i}$, where $\epsilon>0$ is a parameter which
will tend to 0. Since the set of discrete finitely supported
measures is dense in the space of compactly supported Borel
measures on $\Real^n$ in the $w^*$-topology, we may choose
$\mu_\epsilon$ so that as linear functionals, the values of $\mu$
and $\mu_\epsilon$ on the following $n(n+1)/2 + 1$ continuous
functions are $\epsilon$ close:
\[
\abs{\int x_i x_j d\mu_\epsilon(x) - \delta_{i,j}} = \abs{\int x_i
x_j d\mu_\epsilon(x)  - \int x_i x_j d\mu(x)} < \epsilon, \] for
all $1\leq i \leq j \leq n$ and:
\begin{equation} \label{eq:emergency}
\abs{\int \norm{x}^2_K d\mu_\epsilon(x)  - \int \norm{x}^2_K
d\mu(x)} < \epsilon.
\end{equation}
We see that $\mu_\epsilon$ is chosen to be almost isotropic, but
we do not know how to guarantee this in general. Now, repeating
the proof of Lemma \ref{lem:isotropic_is_gaussian}, we see that
$\Lambda_{\mu_\epsilon}$ in (\ref{eq:Lambda}) is a Gaussian r.v.
whose correlation matrix is almost the identity (up to an
$l_\infty$ error of $\epsilon$ w.r.t. the standard basis).
Therefore sending $\epsilon$ to 0, $\Lambda_{\mu_\epsilon}$ tends
to an $n$-dimensional standard Gaussian r.v. almost surely,
implying that $\int \snorm{\sum_{i=1}^{m_\epsilon} g_i(\omega)
\sqrt{\lambda^\epsilon_i} v^\epsilon_i}^2_K d\omega$ tends to
$\int \snorm{\sum_{i=1}^n g_i'(\omega) e_i}^2_K d\omega = n
M_2(K)^2$. Since by the discrete case:
\[
\int \snorm{\sum_{i=1}^{m_\epsilon} g_i(\omega)
\sqrt{\lambda^\epsilon_i} v^\epsilon_i}^2_K d\omega \leq
T_2(X_K)^2 \int\norm{x}^2_K d\mu_\epsilon(x),
\]
and $\int\norm{x}^2_K d\mu_\epsilon(x)$ tends to $\int\norm{x}^2_K
d\mu(x)$ by (\ref{eq:emergency}), this completes the proof.
\end{proof}

One of the most useful isotropic measures associated to the
geometry of a convex body $K$, comes from John's decomposition of
the identity, when $K$ is put in John's maximal volume ellipsoid
position: if $D_n$ is the ellipsoid of maximal volume inside $K$,
there exist contact points $\set{v_i}$ of $D_n$ and $K$ and
positive scalars $\set{\lambda_i}$, such that $\mu_K =
\sum_{i=1}^m \lambda_i \delta_{v_i}$ is isotropic. Since
$\abs{v_i}=1$, it immediately follows that $\sum_{i=1}^m
\lambda_i = n$. Applying Proposition \ref{prop:isotropic-prop}
with the measure $\mu_K$, first with $K$ and then with $K^\circ$,
we immediately have as a corollary the following two known
inequalities. The first essentially appears in
\cite{Maurey-Report}, and in \cite{Milman-Schechtman-Book} with a
worse constant, and the second appears in
\cite{Davis-etal-Lemma}. Both in \cite{Davis-etal-Lemma} and in
\cite{Maurey-Report}, the proofs rely on Operator Theory, whereas
in our approach the elementary geometric flavor is retained, and
both proofs are unified into a single framework.

\begin{cor} \label{cor:2-inequalities}
Let $K$ be a convex body in John's maximal volume ellipsoid
position. Then:
\begin{eqnarray}
\nonumber M_2(K) / b(K) &\geq& 1 / C_2(X_K), \\
\nonumber M^*_2(K) b(K) &\leq& T_2(X^*_K).
\end{eqnarray}
\end{cor}
\begin{proof}
The $b(K)$ terms are simply normalizations to the case that $D_n$
is indeed the ellipsoid of maximal volume inside $K$. It remains
to notice that $\abs{v_i} = \norm{v_i}_K = \norm{v_i}^*_K = 1$,
as contact points between $D_n$ and $K$. Since $\sum_{i=1}^m
\lambda_i = n$, we have:
\[
\brac{\sum_{i=1}^m \lambda_i (\norm{v_i}_K)^2}^{1/2} =
\brac{\sum_{i=1}^m \lambda_i (\norm{v_i}^*_K)^2}^{1/2} = \sqrt{n}.
\]
The assertions now clearly follow from Proposition
\ref{prop:isotropic-prop}.
\end{proof}
\begin{rem}
The other two inequalities:
\begin{eqnarray}
\nonumber M_2(K) / b(K) &\leq& T_2(X_K), \\
\nonumber M^*_2(K) b(K) &\geq& 1/C_2(X^*_K),
\end{eqnarray}
are trivial and loose. The first follows from $M_2(K) \leq b(K)$,
and the second from Urysohn's inequality: \[M^*_2(K) \geq M^*(K)
\geq \VolRatio{K}{D_n} \geq \frac{1}{b(K)}.\]
\end{rem}

By duality, we have:
\begin{cor} \label{cor:2-inequalities-dual}
Let $K$ be a convex body in Lowner's minimal volume outer
ellipsoid position. Then:
\begin{eqnarray}
\nonumber M^*_2(K) / a(K) &\geq& 1 / C_2(X^*_K), \\
\nonumber M_2(K) a(K) &\leq& T_2(X_K).
\end{eqnarray}
\end{cor}

\medskip

Corollary \ref{cor:2-inequalities-dual} shows that having type-2
implies having finite outer volume-ratio (this will be evident in
the proof of the next Theorem), so it is not surprising that we
get the following useful bound on the isotropic constant, when
placing the body in Lowner's outer ellipsoid position. What is a
little more surprising, is that we manage to get the same bound by
putting the body in the isotropic position, and directly applying
Proposition \ref{prop:isotropic-prop} on the (properly
normalized) uniform measure on $K$. The latter part may also be
shown to follow from Theorem 1.4 in \cite{BMMP}.


\begin{thm} \label{thm:L_K-leq-type-2}
Let $K$ be a convex body. Then:
\begin{equation} \label{eq:L_K-leq-type-2}
L_K \leq C \inf \set{\left . T_2(X_L) \VolRatio{L}{K} \right | K
\subset L \textrm{ is a convex body} }.
\end{equation}
In addition, if $\Vol{K} = 1$ and $K$ is in Lowner's minimal
volume outer ellipsoid position or in isotropic position, then:
\[
L_K \leq C \frac{T_2(X_K)}{\sqrt{n} M_2(K)}.
\]
\end{thm}
\begin{proof}
Since (\ref{eq:L_K-leq-type-2}) is invariant under homothety, we
may assume that $\Vol{K}=1$. Now let $L$ be any convex body
containing $K$, and assume that $T(L)$ is in Lowner's minimal
volume outer ellipsoid position, where $T \in SL(n)$. By Corollary
\ref{cor:2-inequalities-dual} and Jensen's inequality (as in
(\ref{eq:Jensen})):
\[
a(T(L)) \leq \frac{T_2(X_L)}{M_2(T(L))} \leq C \sqrt{n}
\VolRatio{L}{K} T_2(X_L).
\]
Using the characterization of $L_K$ mentioned in the Introduction,
we immediately have:
\[
L_K^2 \leq \frac{1}{n} \int_{T(K)} \abs{x}^2 dx \leq \frac{1}{n}
a(T(L))^2 \leq \brac{C \VolRatio{L}{K} T_2(X_L)}^2.
\]

Evidently, the above argument also proves the second part of the
Theorem when $K$ is in Lowner's minimal volume outer ellipsoid
position. When $K$ is in isotropic position, we apply Proposition
\ref{prop:isotropic-prop} to the isotropic measure $d\mu = 1/L_K^2
\chi_K dx$, yielding:
\[
\sqrt{n} M_2(K) \leq T_2(X_K) \frac{1}{L_K} \brac{\int_K
\norm{x}^2_K}^{1/2} \leq T_2(X_K) / L_K.
\]
The assertion therefore follows (even without a constant).
\end{proof}
\begin{rem}
For completeness, it is worthwhile to mention that a different
form of Theorem \ref{thm:L_K-leq-type-2} may be derived from a
deeper result of Milman and Pisier, who showed in
\cite{Milman-Pisier} that the volume-ratio of $K$ is bounded from
above by $C C_2(X_K) \log C_2(X_K)$ (this is an improvement over
the initial bound showed in
\cite{Bourgain-Milman-vr-and-reverse-santalo}). Using another
deep result, the reverse Blaschke-Santalo inequality
(\cite{Bourgain-Milman-vr-and-reverse-santalo}, see
(\ref{eq:reverse-Santalo})), this implies that the outer
volume-ratio of $K$ is bounded from above by $C' C_2(X^*_K) \log
C_2(X^*_K)$, so the same argument as above gives:
\[
L_K \leq C \inf \set{\left . C_2(X^*_L) \log C_2(X^*_L)
\VolRatio{L}{K} \right | K \subset L \textrm{ is a convex body} }.
\]
Since $ C_2(X^*_L) \leq T_2(X_L) \leq C_2(X^*_L) \norm{Rad(X_L)}
$, where $Rad$ denotes the Rademacher projection (see
\cite{Milman-Schechtman-Book}), we see that the two forms are very
similar, but elementary examples show that neither form
out-performs the other.
\end{rem}

\medskip

Since it is well known (e.g. \cite{Milman-Schechtman-Book}) that
subspaces of $L_p$, for $p \geq 2$, have a type-2 constant of the
order of $\sqrt{p}$ (this is a consequence of Kahane's
inequality), we immediately have the following Corollary of
Theorem \ref{thm:L_K-leq-type-2}.

\begin{cor} \label{cor:L_K-leq-sqrt-p-geq-2}
\[
L_K \leq C \inf \set{\left . \sqrt{p} \VolRatio{L}{K} \right | K
\subset L \; , \; L \in SL_p^n \; , \; p \geq 2}.
\]
\end{cor}

\medskip

We conclude this section by giving another application of
Proposition \ref{prop:isotropic-prop}. In principle, it seems
useful to apply it on any isotropic measure which is naturally
associated to a convex body in certain special positions.
Fortunately, in \cite{Giannopoulos-Milman}, Giannopoulos and
Milman have derived a framework to generate such measures, by
considering bodies in minimum quermassintegral positions. We will
only give the following application for the minimal surface-area
position, i.e. the position for which $\Vol{\partial T(K)}$ is
minimal for all $T\in SL(n)$, which was characterized by Petty in
\cite{Petty}. Recall that $\sigma_K$, the area measure of $K$ is
defined on $S^{n-1}$ as:
\[
\sigma_K(A) = \nu\brac{\set{x \in \partial K | n_K(x) \in A}},
\]
where $n_K(x)$ denotes an outer normal to $K$ at $x$ and $\nu$ is
the $n-1$ dimensional surface measure on $K$.


\begin{prop}
Let $K$ be a convex body in minimal surface-area position. Then:
\[
1/C_2(X_K) \leq \frac{M_2(K)}{ \brac{1 / \Vol{\partial K}
\int_{S^{n-1}} \norm{x}_K^2 d\sigma_K(x)}^{1/2}}  \leq T_2(X_K).
\]
\end{prop}
\begin{proof}
It was shown in \cite{Petty} that $K$ is in minimal surface-area
position iff $n / \Vol{\partial K} d\sigma_K$ is isotropic.
Applying Proposition \ref{prop:isotropic-prop} with $\sigma_K$
yields the claimed inequalities.
\end{proof}




\section{Sections and Quotients of $L_p$} \label{sec:SectionsOfLp}

As seen in the previous section, it is actually pretty
straightforward to obtain a bound on the isotropic constant of
any convex body $K$ for which we have control over $T_2(X_K)$,
since in that case $K$ has bounded outer volume-ratio. In
particular, this applies for sections of $L_p$, at least for $p
\geq 2$. In this section, we introduce a new technique involving
dual mixed-volumes, which is well adapted to deal specifically
with integral representations of $\norm{\cdot}^t$. This is well
suited for dealing with sections of $L_p$, since by a classical
result of P. L\'{e}vy (\cite{Levy}), $L \in SL_p^n$ for $p\geq 1$
iff there exists a non-negative Borel measure $\mu_L$ on $S^{n-1}$
such that:
\begin{equation} \label{eq:L_p-representation}
\norm{x}_L^p = \int_{S^{n-1}} \abscalar{x,\theta} ^p
d\mu_L(\theta),
\end{equation}
for all $x \in \Real^n$. This characterization extends to any
$p>0$, and it will enable us to extend the bound on $L_K$ to the
case $K \in SL_p^n$ for all $p > 0$. As we shall see, for a
general convex body $K$, it is not the \emph{volume-ratio}
between $L \in SL_p^n$ containing $K$ and $K$ which matters, but
rather some other natural parameter. Moreover, our new technique
will enable us to pass to the dual, and recover Junge's bound on
the isotropic constant of quotients of $L_q$. In Section
\ref{sec:k-BP-bodies}, we continue to apply our technique to
bound the isotropic constant of convex bodies contained in
$k$-Busemann-Petty bodies.

\begin{thm} \label{thm:main1-basic}
Let $K$ be a centrally symmetric convex body in isotropic
position, and let $D$ be a Euclidean ball normalized so that
$\Vol{D} = \Vol{K}$. Then for any $p > 0$ and any $L \in SL_p^n$:
\begin{equation} \label{eq:main1-basic}
C_1 / \sqrt{p_0} \; \leq \; L_K \; / \brac{ \frac{
\widetilde{V}_{-p}(L,K) } {\widetilde{V}_{-p}(L,D)}}^{1/p} \leq
C_2 \sqrt{p_0},
\end{equation}
where $p_0 = \max(1,\min(p,n))$.
\end{thm}

\begin{rem}
By taking the limit in (\ref{eq:L_p-representation}) as
$p\rightarrow 0+$, we may define $SL_0^n$ to be the class of
$n$-dimensional star-bodies $L$ for which:
\[
\norm{x}_L = \exp\brac{\int_{S^{n-1}} \log \abscalar{x,\theta}
d\mu_L(\theta) + C},
\]
for some Borel probability measure $\mu_L$ and constant $C$ and
all $x\in \Real^n$. In that case, Theorem \ref{thm:main1-basic}
holds true for $p=0$ as well (by passing to the limit), if we
replace the expressions of the form
$\widetilde{V}_{-p}(L_1,L_2)^{1/p}$ appearing in (4.2), by the
limit as $p \rightarrow 0+$ assuming $\Vol{L_2} = 1$, namely
$\exp\brac{1/n \int_{S^{n-1}} \log(\rho_{L_2}(x) / \rho_{L_1}(x))
\rho_{L_2}(x)^n dx}$.
\end{rem}


\begin{proof}[Proof of Theorem \ref{thm:main1-basic}]

Let $\mu_L$ denote the Borel measure on $S^{n-1}$ from
(\ref{eq:L_p-representation}) corresponding to $L$. Then for any
star-body $G$:
\begin{eqnarray}
\nonumber \widetilde{V}_{-p}(L,G) & = & \frac{1}{n} \int_{S^{n-1}}
\norm{x}_L^p \norm{x}_G^{-(n+p)} dx \\
\nonumber & = & \frac{1}{n} \int_{S^{n-1}} \int_{S^{n-1}} \abscalar{x,\theta}^p d\mu_L(\theta) \norm{x}_G^{-(n+p)} dx \\
\nonumber & = & \frac{1}{n} \int_{S^{n-1}} d\mu_L(\theta)
\int_{S^{n-1}} \abscalar{x,\theta}^p \norm{x}_G^{-(n+p)} dx \\
\label{eq:fubini1} & = & \frac{n+p}{n} \int_{S^{n-1}}
d\mu_L(\theta) \int_{G} \abscalar{x,\theta}^p dx
\end{eqnarray}

Let us evaluate the expression $\int_{G} \abscalar{x,\theta}^p
dx$. If $G$ is of volume 1 and $p\geq 1$, then by Jensen's
inequality:
\begin{equation} \label{eq:lp-jensen}
\int_G \abscalar{x,\theta} dx \leq \brac{ \int_G
\abscalar{x,\theta}^p dx }^{1/p} \;\;\; \forall p\geq 1.
\end{equation}
If $G$ is in addition convex, then by a well known consequence of
a lemma by C. Borell (\cite{Borell-lemma}), it follows that the
linear functional $\scalar{\cdot,\theta}$ has a $\psi_1$-type
behaviour on $G$, and therefore:
\begin{equation} \label{eq:psi1}
\brac{ \int_G \abscalar{x,\theta}^p dx } ^ {1/p} \leq C p \int_G
\abscalar{x,\theta} dx \;\;\; \forall p \geq 1
\end{equation}
If in addition $p\geq n$, it is well known that (e.g. \cite[Lemma
4.1]{Paouris-Small-Diameter}):
\begin{equation} \label{eq:max-norm}
\brac{ \int_G \abscalar{x,\theta}^p dx } ^ {1/p} \simeq
\norm{\theta}^*_G \;\;\; \forall p \geq n.
\end{equation}
Finally, if $G$ is convex, of volume 1 and $0<p<1$, then it
follows from the estimates in Corollary 2.5 and 2.7 in
\cite{Milman-Pajor-LK} that:
\begin{equation} \label{eq:p-leq-1}
\brac{ \int_G \abscalar{x,\theta}^p dx } ^ {1/p} \simeq \int_G
\abscalar{x,\theta} dx \;\;\; \forall p\in(0,1).
\end{equation}

The expression in (\ref{eq:main1-basic}) is invariant under
simultaneous homothety of $K$ and $D$, so we may assume that
$\Vol{K} = \Vol{D} = 1$. Since $K$ is in isotropic position, we
have $\int_K \scalar{x,\theta}^2 dx = L_K^2$ for all $\theta \in
S^{n-1}$, and by (\ref{eq:lp-jensen}) - (\ref{eq:p-leq-1}) it
follows that for all $\theta \in S^{n-1}$:
\begin{equation} \label{eq:isotropic-psi1}
A \leq \brac{ \int_K \abscalar{x,\theta}^p dx } ^ {1/p} / \; L_K
\leq B p_0 \; \; \; \forall p > 0.
\end{equation}
It remains to notice that for a Euclidean ball $D$ of volume 1, a
straightforward computation (in the case $1\leq p \leq n$)
together with (\ref{eq:max-norm}) and (\ref{eq:p-leq-1}), gives
that for all $\theta \in S^{n-1}$:
\begin{equation} \label{eq:euclidean-ball}
\brac{ \int_D \abscalar{x,\theta}^p dx } ^ {1/p} \simeq \sqrt{p_0}
\; \; \; \forall p>0.
\end{equation}
By (\ref{eq:fubini1}), we have:
\[
\brac{ \frac{ \widetilde{V}_{-p}(L,K) }
{\widetilde{V}_{-p}(L,D)}}^{1/p} = \brac{\frac{\int_{S^{n-1}}
d\mu_L(\theta) \int_{K} \abscalar{x,\theta}^p dx}{\int_{S^{n-1}}
d\mu_L(\theta) \int_{D} \abscalar{x,\theta}^p dx}}^{1/p}.
\]
Since $\mu_L \geq 0$, using (\ref{eq:isotropic-psi1}) and
(\ref{eq:euclidean-ball}), we get the required
(\ref{eq:main1-basic}):
\[
\frac{1}{C_2 \sqrt{p_0}} \leq \brac{ \frac{
\widetilde{V}_{-p}(L,K) } {\widetilde{V}_{-p}(L,D)}}^{1/p} / L_K
\leq \frac{\sqrt{p_0}}{C_1}.
\]
\end{proof}

\begin{rem}
Notice that for $0 \leq p < 1$, the unit-ball of a subspace of
$L_p$ is no longer necessarily a convex body. We will see more
examples where $L$ is a non-convex star-body later on. In fact,
using the results in \cite{Guedon-extension-to-negative-p} of
Gu\'{e}don, it is possible to extend Theorem
\ref{thm:main1-basic} to $p>-1$, but then the constants $C_1$ and
$C_2$ will depend on $p$. We do not proceed in this direction,
because we are able to show in Section \ref{sec:k-BP-bodies} that
Theorem \ref{thm:main1-basic} is also valid for $p=-1$ (then
$SL_p^n$ is replaced by the class of intersection-bodies), and we
are able to generalize this to $k$-Busemann-Petty bodies.
\end{rem}

\medskip

We can now extend Corollary \ref{cor:L_K-leq-sqrt-p-geq-2} to the
following more general result.

\begin{thm} \label{thm:main1-particular}
Let $K$ be a centrally symmetric convex body in isotropic
position with $\Vol{K}=\Vol{D_n}$. Then:
\[
L_K \leq C \inf \set{\left . \frac{\sqrt{p_0}}{M_p(L)}  \right | K
\subset L \; , \; L \in SL_p^n \; , \; p \geq 0},
\]
where $p_0 = \max(1,\min(p,n))$.
\end{thm}

\begin{proof}
If $K \subset L$, then obviously
$\widetilde{V}_{-p}(L,K) \leq \widetilde{V}_{-p}(K,K) = \Vol{K}$.
Applying Theorem \ref{thm:main1-basic} with $\Vol{D} = \Vol{K} =
\Vol{D_n}$, (\ref{eq:main1-basic}) implies:
\begin{equation} \label{eq:main1-particular}
L_K \leq C_2 \sqrt{p_0} \brac{\frac{\Vol{D_n}}{\frac{1}{n}
\int_{S^{n-1}} \rho_L(x)^{-p} dx}}^{1/p} = C_2
\frac{\sqrt{p_0}}{M_p(L)}.
\end{equation}
\end{proof}

Using Jensen's inequality (\ref{eq:Jensen}) and homogeneity, we
immediately have the following corollary, which unifies the
bounds on $L_K$ for $SL_p^n$ of Ball (the case $1 \leq p \leq 2$)
and Junge (the case $p\geq 2$), and extends their results to $p
\geq 0$:

\begin{cor} \label{cor:main1-particular}
For any centrally symmetric convex body $K$:
\[
L_K \leq C \inf \set{\left . \sqrt{p_0} \VolRatio{L}{K} \right | K
\subset L \; , \; L \in SL_p^n \; , \; p \geq 0},
\]
where $p_0 = \max(1,\min(p,n))$.
\end{cor}


\begin{rem} \label{rem:psi_2}
Notice that the proof of Theorem \ref{thm:main1-basic} does not
use the assumption that the body $D$ is a Euclidean ball: the only
property used is the one in (\ref{eq:euclidean-ball}). In fact,
for the right-hand inequality in (\ref{eq:main1-basic}), $D$ may
be chosen as any $\psi_2$-body in isotropic position. Recall that
$D$ is called a $\psi_2$-body (with constant $A>1$), if for all $p
\geq 1$:
\[
\brac{\int_D \abs{\scalar{x,\theta}}^p dx}^{1/p} \leq A \sqrt{p}
\brac{\int_D \abs{\scalar{x,\theta}}^2 dx}^{1/2} \;\; \; \forall
\theta \in S^{n-1}.
\]
Bourgain has shown in \cite{Bourgain-Psi-2-Bodies} that if $D$ is
a $\psi_2$-body then $L_D \leq C A \log{A}$. Therefore if $D$ is a
$\psi_2$-body of volume 1 in isotropic position,
(\ref{eq:euclidean-ball}) may be replaced by:
\[
\brac{\int_D \abs{\scalar{x,\theta}}^p dx}^{1/p} \leq A^2 \log{A}
\sqrt{p}  \;\;\; \; \forall \theta \in S^{n-1} \; , \; \forall p
\geq 1.
\]
(\ref{eq:main1-particular}) then reads (when
$\Vol{K}=\Vol{D}=\Vol{D_n}$):
\[
L_K \leq C(A) \sqrt{p_0}
\brac{\frac{\Vol{D_n}}{\widetilde{V}_{-p}(L,D)}}^{1/p} = C(A)
\sqrt{p_0} \brac{ \int_{S^{n-1}} \norm{x}^p_L \rho_D(x)^{n+p}
d\sigma(x)}^{-1/p}.
\]

By Bourgain's result, if all linear functionals are $\Psi_2$,
then $L_K$ is bounded. Ironically, it follows from the proof of
Theorem \ref{thm:main1-basic} that if all linear functionals have
``bad" $\psi_2$ behaviour, e.g.
\[
\brac{\int_K \abs{\scalar{x,\theta}}^q dx}^{1/q} \geq C \sqrt{q}
\int_K \abs{\scalar{x,\theta}} dx \;\; \forall \theta \in S^{n-1},
\]
for a certain $q\geq 1$, then the bound on $L_K$ improves ($L_K
\leq C \VolRatio{L}{K}$ for all $L\in SL_q^n$ containing $K$, in
the example above). Perhaps this may be used to our advantage?
\end{rem}


\medskip

We now turn to reproduce Junge's bound on $L_K$ for quotients of
$L_q$. As mentioned in the Introduction, for $1 < q \leq 2$,
Junge's result is more general than ours and applies to all
\emph{subspaces} of quotients of $L_q$. Nevertheless, our proof
provides a (formally) stronger bound, applies to the entire range
$1 < q \leq \infty$, and retains the problem's Geometric nature,
avoiding unnecessary tools from Operator Theory.
In addition, although derived independently, our proof is very
similar to Bourgain's proof that $L_K \leq C n^{1/4} \log(1+n)$,
and the latter may be thought of as an extremal case of our
proof, where our argument breaks down.

\begin{thm} \label{thm:main1-dual}
Let $K$ be a centrally symmetric convex body in isotropic
position with $\Vol{K}=\Vol{D_n}$. Then:
\[
L_K \leq C \; \inf \set{ \sqrt{p_0} \; M^*_p(T(L)) \; \left | \:
\begin{array}{c} K \subset L \; , \; L \in QL_q^n \; , \; T \in
SL(n) \; , \\  1\leq q\leq\infty \; , \; 1/p + 1/q = 1
\end{array} \right . },
\]
where $p_0 = \min(p,n)$ and $p = q^*$ is the conjugate exponent
to $q$.
\end{thm}

We postpone the proof of Theorem \ref{thm:main1-dual} for later.
In order to see why this Theorem implies Junge's bound for
quotients of $L_q$, we will need the following lemma:




\begin{lem} \label{lem:M_p-bound}
Let $K$ be a convex body with $\Vol{K} = \Vol{D_n}$.
\begin{enumerate}
\item
If $K \in SL_p^n$ for $1\leq p \leq \infty$, then there exists a
position of $K$ for which $M_p(K) \leq C \sqrt{p_0}$, where $p_0 =
\min(p,n)$.
\item
If $K \in QL_q^n$ for $1 \leq q \leq \infty$, then there exists a
position of $K$ for which $M^*_p(K) \leq C \sqrt{p_0}$, for $p =
q^* = q/(q-1)$ and $p_0$ as above.
\end{enumerate}
\end{lem}

Applying the second part of the lemma to the body $L$ from
Theorem \ref{thm:main1-dual} and using homogeneity, we immediately
have:

\begin{cor} \label{cor:main1-dual}
For any centrally symmetric convex body $K$:
\[
L_K \leq C \inf \set{\left . p_0 \VolRatio{L}{K} \: \right |
\begin{array}{c} K \subset L \; , \; L \in QL_q^n \; , \\ 1 \leq q \leq
\infty \; , \; 1/p + 1/q = 1 \end{array} },
\]
where $p_0 = \min(p,n)$.
\end{cor}

\begin{proof}[Proof of Lemma \ref{lem:M_p-bound}]
We will prove part 1. Part 2 then follows easily by duality,
using the reverse Blaschke-Santalo inequality
(\cite{Bourgain-Milman-vr-and-reverse-santalo}):
\begin{equation} \label{eq:reverse-Santalo}
\VolRatio{K}{D_n} \VolRatio{K^\circ}{D_n} \geq c,
\end{equation}
to ensure that the volume of $K^\circ$ is not too small.

The case $1\leq p \leq 2$ is straightforward, since for this range
it is well known that sections of $L_p$ have finite volume-ratio
(for instance, because they have cotype-2 and using
\cite{Bourgain-Milman-vr-and-reverse-santalo}, or by
\cite{Ball-Volume-Ratios}). Therefore, in John's maximal volume
ellipsoid position, $M_p(K) \leq b(K) \leq C$. We remark that it
remains to prove the lemma for $2 \leq p \leq n$, since it is
known that $M_p(K) \simeq M_n(K) \simeq b(K)$ for $p>n$ (e.g.
\cite{LMS}).

We will present three different proofs for the case $2\leq p \leq
n$, placing the body $K$ in three different positions. We note
that the first two proofs actually prove a stronger statement:
for any $K \in SL_p^n$ there exists a position in which $M_p(K)
\leq C \sqrt{p} / a(K)$. Since this formulation is volume free,
we do not really need the reverse Blaschke-Santalo inequality to
prove the dual second part of the lemma (for the range $1\leq
q\leq 2$). The third proof is an elementary consequence of Theorem
\ref{thm:main1-particular}, and appears also in Corollary
\ref{cor:bounds-on-mean-radius}.


\begin{enumerate}
\item
If $2 \leq p \leq n$, then $T_2(X_K) \leq C \sqrt{p}$ (by Kahane's
inequality), so by Corollary \ref{cor:2-inequalities-dual}, if
$K$ is in \emph{Lowner's minimal volume outer ellipsoid position},
then $M_2(K) a(K) \leq C \sqrt{p}$. Notice that in Lowner's
position, $b(K) \leq \sqrt{n} / a(K)$. Since $\Vol{K} =
\Vol{D_n}$, we obviously have $a(K) \geq 1$, implying that
$M_2(K) \leq C \sqrt{p}$ and $b(K) \leq \sqrt{n}$. We now use a
known result from \cite{LMS}, stating that $M_p(K) \simeq
\max(M_2(K),b(K)\sqrt{p}/\sqrt{n})$, which under our conditions
implies $M_p(K) \leq C \sqrt{p}$.
\item
By approximation, we may assume that $K$ is a section of $l_p^m$,
for some large enough $m$. We will put $K$ in the \emph{Lewis
position} (\cite{Lewis-Position}), as used in
\cite{Ball-Volume-Ratios}. In this position, there exists a
sequence of $m$ unit vectors $\set{u_i}$ and positive scalars
$\set{c_i}$, such that $\norm{x}^p_K = \sum_{i=1}^m c_i
\abs{\scalar{x,u_i}}^p$ and such that $\mu = \sum_{i=1}^m c_i
\delta_{u_i}$ is an isotropic measure (see Section
\ref{sec:type-2}). In particular, $\sum_{i=1}^m c_i = n$. An
elementary computation shows that for $2\leq p \leq n$:
\[
M_p(K) = \brac{\sum_{i=1}^m c_i \int_{S^{n-1}}
\abs{\scalar{\theta,u_1}}^p d\sigma(\theta)}^{1/p}  \simeq
\brac{\sum_{i=1}^m c_i}^{1/p} \frac{\sqrt{p}}{\sqrt{n}} =
\frac{\sqrt{p}}{n^{1/2 - 1/p}}.
\]
But in this position, H\"{o}lder's inequality shows that:
\[
\abs{x}^2 = \sum_{i=1}^m c_i \abs{\scalar{x,u_i}}^2 \leq
\brac{\sum_{i=1}^m c_i}^{1-2/p} \brac{\sum_{i=1}^m c_i
\abs{\scalar{x,u_i}}^p}^{2/p} = n^{1-2/p} \norm{x}^2_K,
\]
and therefore $a(K) \leq n^{1/2 - 1/p}$. It follows that $M_p(K)
\leq C \sqrt{p} / a(K)$, as required.
\item
Put the body $K$ in \emph{isotropic position}, and apply Theorem
\ref{thm:main1-particular} with $L=K$. Using the well known fact
that $L_K$ is always bounded from below by a universal constant
(e.g. \cite{Milman-Pajor-LK}), we immediately have $M_p(K) \leq C
\sqrt{p_0} \sVolRatio{K}{D_n}$, and this is valid for all $p \geq
0$, with $p_0 = \max(1,\min(p,n))$.
\end{enumerate}


\end{proof}

\begin{proof}[Proof of Theorem \ref{thm:main1-dual}]
Let $K$ be in isotropic position and assume $\Vol{K}=1$. Fix $q>1$
and let $L \in QL_q^n$ contain $K$. By duality, $L^\circ$, the
polar body to $L$, is a section of $L_p$, and so is $T(L^\circ)$
for any $T \in SL(n)$. Applying Theorem \ref{thm:main1-basic}, the
left (!) hand side of (\ref{eq:main1-basic}) gives:
\begin{equation} \label{eq:main1-polar}
L_K  \sqrt{p_0} / C_1 \geq \brac{ \frac{
\widetilde{V}_{-p}(T(L^\circ),K) }
{\widetilde{V}_{-p}(T(L^\circ),D)}}^{1/p} \geq \brac{ \frac{
\widetilde{V}_{-p}(T(K^\circ),K) }
{\widetilde{V}_{-p}(T(L^\circ),D)}}^{1/p},
\end{equation}
for $D$ the Euclidean ball of volume 1. Evaluating the numerator
on the right using the trivial $\norm{x}_{T(K^\circ)} \norm{x}_K
\geq \abs{\scalar{T^{-1}(x),x}}$, we have that for any
positive-definite $T \in SL(n)$:
\begin{eqnarray}
\nonumber & & \brac{ \widetilde{V}_{-p}(T(K^\circ),K) }^{1/p} =
\brac{\frac{1}{n} \int_{S^{n-1}} \norm{x}^p_{T(K^\circ)}
\norm{x}_K^{-(n+p)} dx}^{1/p} \\
\nonumber & \geq & \brac{\frac{1}{n} \int_{S^{n-1}}
\abs{\scalar{T^{-1}(x),x}}^p \norm{x}_K^{-(n+2p)}dx}^{1/p} \\
\nonumber & = & \brac{\frac{n+2p}{n} \int_K
\abs{\scalar{T^{-1}(x),x}}^p dx}^{1/p} \geq \int_K
\scalar{T^{-1}(x),x} dx \\
\nonumber & = & tr(T^{-1}) L_K^2 \geq det(T^{-1})^{1/n} n L_K^2 =
n L_K^2,
\end{eqnarray}
where we have used Jensen's inequality, the fact that $\int_K x_i
x_j dx = L_K^2 \delta_{i,j}$, and the Arithmetic-Geometric means
inequality (since $T$ is positive-definite). Together with
(\ref{eq:main1-polar}), and cancelling out one $L_K$ term, this
gives:
\begin{eqnarray}
\nonumber L_K & \leq & \frac{\sqrt{p_0}}{C_1 n}
\brac{\widetilde{V}_{-p}(T(L^\circ),D)}^{1/p} \\
\nonumber & = & \frac{\sqrt{p_0}}{C_1 n} \Vol{D_n}^{-1/n}
M_p(T(L^\circ)) \simeq \frac{\sqrt{p_0}}{\sqrt{n}}
M^*_p((T^{-1})^*(L)),
\end{eqnarray}
for any $T \in SL(n)$ (since it can be factorized into a
composition of a rotation and a positive-definite transformation,
and $M_p$ is invariant to rotations). Changing normalization from
$\Vol{K}=1$ to $\Vol{K}=\Vol{D_n}$, we have the desired:
\[
L_K \leq C \sqrt{p_0} M^*_p(T(L)).
\]
\end{proof}

\begin{rem}
As already mentioned, the proof of Theorem \ref{thm:main1-dual}
clearly resembles Bourgain's proof that $L_K \leq C n^{1/4}
\log(1+n)$. In this respect, we mention that instead of using
$\sqrt{p_0}$ on the left hand side of (\ref{eq:main1-basic}) or
$\sqrt{p}$ on the left hand side of (\ref{eq:main1-polar}), it is
easy to check that one may use $A$, if $K$ is a $\Psi_2$ body with
constant $A$ (as defined in Remark \ref{rem:psi_2}). This implies
that whenever $A < \sqrt{p}$, we get a better bound on $L_K$.
Bourgain has shown that in the general case, one may always
assume that $A \leq n^{1/4}$, but this does not seem to help us in
our context.
\end{rem}

\medskip

To conclude this section, we mention that for a \emph{general}
convex body $K$ (not necessarily a section of $L_p$),
representations other than (\ref{eq:L_p-representation}) of
$\norm{\cdot}_K$ as a spherical convolution of a kernel with a
non-negative Borel measure on $S^{n-1}$ are known. Repeating the
relevant parts of the proof of Theorem \ref{thm:main1-basic} with
$L=K$, it may be possible to bound some natural parameter of the
body $K$ other than $L_K$.



\section{$k$-Busemann-Petty bodies} \label{sec:k-BP-bodies}

An analogous result to Theorem \ref{thm:main1-basic} for
$k$-Busemann-Petty bodies is the following:

\begin{thm} \label{thm:main2-basic}
Let $K$ be a centrally symmetric convex body in isotropic
position, and let $D$ be a Euclidean ball normalized so that
$\Vol{D} = \Vol{K}$. Then for any integer $k = 1,\ldots,n-1$ and
any $L \in \BP_k^n$:
\begin{equation} \label{eq:main2-basic}
C_1 \; \leq \; L_K \; / \brac{ \frac{ \widetilde{V}_{k}(L,D) }
{\widetilde{V}_{k}(L,K)}}^{1/k} \leq C_2 \Lk.
\end{equation}
\end{thm}


\begin{proof}
By definition, if $L \in \BP_k^n$ there exists a Borel measure
$\mu_L$ on $G(n,n-k)$ such that:
\begin{equation} \label{eq:BP_k^n-representation}
\rho_L^k = R_{n-k}^{*} (d\mu_L).
\end{equation}
Therefore, for any star-body $G$:
\begin{eqnarray}
\nonumber \widetilde{V}_{k}(L,G) & = & \frac{1}{n} \int_{S^{n-1}}
\rho_L(x)^k \rho_G(x)^{n-k} dx \\
\nonumber & = & \Vol{D_n} \int_{S^{n-1}} R_{n-k}^{*} (d\mu_L)
(x) \rho_G(x)^{n-k} d\sigma(x) \\
\nonumber & = & \Vol{D_n} \int_{G(n,n-k)} R_{n-k}
(\rho_G^{n-k}) (E) d\mu_L(E) \\
\label{eq:fubini2} & = & \frac{\Vol{D_n}}{\Vol{D_{n-k}}}
\int_{G(n,n-k)} \Vol{G \cap E} d\mu_L(E).
\end{eqnarray}

The expression in (\ref{eq:main2-basic}) is invariant under
simultaneous homothety of $K$ and $D$, so we may assume that
$\Vol{K} = \Vol{D} = 1$. It is known
(\cite{Ball-PhD},\cite{Milman-Pajor-LK},\cite{Ball-kdim-sections})
that for a convex $K$ in isotropic position and volume 1:
\begin{equation} \label{eq:sections-of-isotropic}
A \leq \Vol{K \cap E}^{1/k} L_K \leq B \Lk \; \; \; \forall E \in
G(n,n-k).
\end{equation}
The proof of (\ref{eq:sections-of-isotropic}) is based on the
fact that the function $f(x) = \Vol{K \cap \set{E + x}}$ on
$E^\perp$ is log-concave and isotropic, and its isotropic constant
is $L_f = f(0)^{1/k} L_K$. It was shown in (\cite{Ball-PhD}) that
an isotropic log-concave function $f$ on $\Real^k$ satisfies $A
\leq L_f \leq B \Lk$, implying (\ref{eq:sections-of-isotropic}).

It remains to notice that for a Euclidean ball $D$ of volume 1, a
straightforward computation shows that for any $k=1,\ldots,n-1$:
\begin{equation} \label{eq:sections-of-ball}
\Vol{D \cap E} ^ {1/k} \simeq 1 \; \; \; \forall E \in G(n,n-k).
\end{equation}
By (\ref{eq:fubini2}), we have:
\[
\brac{ \frac{ \widetilde{V}_{k}(L,K) }
{\widetilde{V}_{k}(L,D)}}^{1/k} = \brac{\frac{\int_{G(n,n-k)}
\Vol{K \cap E} d\mu_L(E)}{\int_{G(n,n-k)} \Vol{D \cap E}
d\mu_L(E)}}^{1/k}.
\]
Since $\mu_L \geq 0$, using (\ref{eq:sections-of-isotropic}) and
(\ref{eq:sections-of-ball}), we get the required
(\ref{eq:main2-basic}):
\[
C_1 \leq \brac{ \frac{ \widetilde{V}_{k}(L,K) }
{\widetilde{V}_{k}(L,D)}}^{1/k}  L_K \leq C_2 \Lk.
\]
\end{proof}

\begin{rem}
It is known (\cite{Koldobsky-I-equal-BP}) that the representation
(\ref{eq:BP_k^n-representation}) exists for any star-body $L$
whose radial function $\rho_L$ is infinitely times differentiable
on $S^{n-1}$, if we allow $\mu_L = \mu_{L,k}$ to be a
\emph{signed} measure on $G(n,n-k)$. Using $L=K$ for example, and
repeating the argument in the proof of Theorem
\ref{thm:main2-basic}, we get that:
\[
L_K \leq C \; \brac{\frac{\int_{G(n,n-k)} \abs{d\mu_{K,k}}(E)}
{\int_{G(n,n-k)} d\mu_{K,k}(E)}}^{1/k} \L_k,
\]
so it remains to evaluate the above ratio. Unfortunately, this
approach does not seem promising, since for a general smooth
function $f$ on $S^{n-1}$, for which the representation $f =
R_{n-k}^*(d\mu)$ is known to exist, it is easy to show that this
ratio may be arbitrarily large for $k=1$ and a fixed value of $n$.
\end{rem}

\medskip

We can now prove analogous results to Theorem
\ref{thm:main1-particular} and Corollary
\ref{cor:main1-particular}.

\begin{thm} \label{thm:main2-particular}
Let $K$ be a centrally symmetric convex body in isotropic
position with $\Vol{K}=\Vol{D_n}$. Then:
\[
L_K \leq C \inf \set{\left . \L_k \MR_k(L)  \right | K \subset L
\; , \; L \in \BP_k^n \; , \; k=1,\ldots,n-1}.
\]
\end{thm}
\begin{proof}
If $K \subset L$, then obviously $\widetilde{V}_{k}(L,K) \geq
\widetilde{V}_{k}(K,K) = \Vol{K}$. Applying Theorem
\ref{thm:main2-basic} with $\Vol{D} = \Vol{K} = \Vol{D_n}$,
(\ref{eq:main2-basic}) implies:
\begin{equation} \label{eq:main2-particular}
L_K \leq C_2 \Lk \brac{\frac{\frac{1}{n} \int_{S^{n-1}}
\rho_L(x)^{k} dx}{\Vol{D_n}}}^{1/k}  = C_2 \Lk \MR_k(L).
\end{equation}
\end{proof}

Using Jensen's inequality (\ref{eq:Jensen}) and homogeneity, we
immediately have the following corollary, which generalizes
Ball's bound on $L_K$ for $SL_p^n$ with $1 \leq p \leq 2$, since
in that range $SL_p^n \subset \BP_k^n$ for \mbox{$k=1,\ldots,n-1$}
(as explained in the Introduction):

\begin{cor} \label{cor:main2-particular}
For any centrally symmetric convex body $K$:
\[
L_K \leq C \inf \set{\left . \L_k \VolRatio{L}{K} \right | K
\subset L \; , \; L \in \BP_k^n \; , \; k=1,\ldots,n-1}.
\]
\end{cor}

\begin{rem}
As before, the proof of Theorem \ref{thm:main2-basic} does not
utilize the assumption that $D$ is a Euclidean ball. The only
property of $D$ used is the one stated in
(\ref{eq:sections-of-ball}). By a result of Junge
(\cite{Junge-Unconditional}), this is satisfied by any
1-unconditional convex body in isotropic position.
(\ref{eq:main2-particular}) then reads (when $\Vol{K} = \Vol{D} =
\Vol{D_n}$):
\[
L_K \leq C_2 \Lk
\brac{\frac{\widetilde{V}_{k}(L,D)}{\Vol{D_n}}}^{1/k} = C_2 \Lk
\brac{\int_{S^{n-1}} \rho_L(x)^k \rho_D(x)^{n-k}
d\sigma(x)}^{1/k}.
\]
\end{rem}

\medskip


As in the previous section, we may prove dual counterparts to
Theorem \ref{thm:main2-particular} and Corollary
\ref{cor:main2-particular}. Before proceeding, we will need the
following useful lemma:

\begin{lem} \label{lem:average-sections-generic}
For any compact set $A \subset \Real^n$ and $m=1,\ldots,n$:
\[
\int_{G(n,m)} \Vol{A \cap E} d\nu(E) \leq \inf_{T \in SL(n)}
\sup_{E \in G(n,m)} \Vol{T(A) \cap E},
\]
where $\nu$ is the Haar probability measure on $G(n,m)$.
\end{lem}
\begin{proof}
Notice that for any compact set $A \subset \Real^n$ and $T \in
SL(n)$:
\[
\Vol{A \cap E} = D_T(E) \Vol{T(A) \cap T(E)},
\]
where the Jacobian $D_T(E)$ does not depend on $A$. Now let $D$
be the Euclidean ball of volume 1, fix $T \in SL(n)$, and denote
$G = G(n,m)$ for short. Denote $M = \sup_{E \in G} \Vol{T(A) \cap
E}$. Then:
\begin{eqnarray}
\nonumber & & \int_{G} \Vol{A \cap E} d\nu(E) = \int_{G}
\Vol{T(A) \cap T(E)} D_T(E) d\nu(E) \\
\nonumber & \leq & M \int_{G} D_T(E) d\nu(E) = M
\frac{\Vol{D_n}^{m/n}}{\Vol{D_m}} \int_{G} \Vol{D \cap T(E)}
D_T(E) d\nu(E) \\
\nonumber & = & M  \frac{\Vol{D_n}^{m/n}}{\Vol{D_m}} \int_{G}
\Vol{T^{-1}(D) \cap E} d\nu(E).
\end{eqnarray}
Now, using polar coordinates, double integration and Jensen's
inequality, we have:
\begin{eqnarray}
\nonumber &  & \!\!\!\!\!\!\! \int_{G} \Vol{T^{-1}(D) \cap E}
d\nu(E) = \Vol{D_m} \int_G \int_{S^{n-1} \cap E}
\norm{\theta}^{-m}_{T^{-1}(D)}
d\sigma_{E}(\theta) d\nu(E) \\
\nonumber & = & \Vol{D_m} \int_{S^{n-1}}
\norm{\theta}^{-m}_{T^{-1}(D)} d\sigma(\theta) \leq \Vol{D_m}
\brac{\int_{S^{n-1}} \norm{\theta}^{-n}_{T^{-1}(D)}
d\sigma(\theta)}^{\frac{m}{n}} \\
\nonumber & = & \Vol{D_m}
\brac{\frac{\Vol{T^{-1}(D)}}{\Vol{D_n}}}^{\frac{m}{n}} =
\frac{\Vol{D_m}}{\Vol{D_n}^{\frac{m}{n}}}.
\end{eqnarray}
We therefore see that for any $T \in SL(n)$:
\[
\int_{G(n,m)} \Vol{A \cap E} d\nu(E) \leq \sup_{E \in G} \Vol{T(A)
\cap E},
\]
which proves the assertion.
\end{proof}

\begin{rem}
An alternative way to prove Lemma
\ref{lem:average-sections-generic} was suggested to us by the
referee, to whom we are grateful. It makes use of a very
interesting result by Grinberg
(\cite{Grinberg-Affine-Invariant}), which was unknown to this
author. In hope of interesting the unfamiliar reader, we bring it
here. The dual affine Quermassintegral of a compact set $A$, which
was introduced by Lutwak in the 80's (see also
\cite{Lutwak-Harmonic}), is defined (up to normalization) as:
\[
\Phi_{n-m}(A) = \brac{\int_{G(n,m)} \Vol{A \cap E}^n
d\nu(E)}^{1/n}.
\]
It was shown in \cite{Grinberg-Affine-Invariant} that
$\Phi_{n-m}$ is indeed invariant to volume preserving linear
transformations: $\Phi_{n-m}(T(A)) = \Phi_{n-m}(A)$ for all $T
\in SL(n)$. Using this, Lemma \ref{lem:average-sections-generic}
is easily deduced from Jensen's inequality, since for any $T\in
SL(n)$:
\begin{eqnarray}
\nonumber \int_{G(n,m)} \Vol{A \cap E} d\nu(E) \leq
\brac{\int_{G(n,m)} \Vol{A \cap E}^n d\nu(E)}^{1/n} \\
\nonumber = \Phi_{n-m}(A) = \Phi_{n-m}(T(A)) \leq \sup_{E \in
G(n,m)} \Vol{T(A) \cap E}.
\end{eqnarray}
We mention another result from \cite{Grinberg-Affine-Invariant},
stating that for a convex body $K$:
\[
\Phi_{n-m}(K) \leq C_{m,n} \Vol{K}^{m/n},
\]
where $C_{m,n}$ is determined by choosing $K = D_n$, and with
equality iff $K$ is a centrally symmetric ellipsoid. This may be
used to give a universal bound for the expression appearing in
the next Lemma \ref{lem:average-sections}, but we will need an
estimate depending on $L_K$ for the proof of Theorem
\ref{thm:main2-dual}.
\end{rem}

Applying Lemma \ref{lem:average-sections-generic} on a convex body
$K$ of volume 1, and using (\ref{eq:sections-of-isotropic}) when
$T(K)$ is in isotropic position, we immediately get the following
lemma as a corollary:

\begin{lem} \label{lem:average-sections}
For any centrally symmetric convex body $K$ with
\mbox{$\Vol{K}=1$:}
\[
\brac{\int_{G(n,n-k)} \Vol{K \cap E} d\nu(E)}^{1/k} \leq C \Lk /
L_K,
\]
where $\nu$ is the Haar probability measure on $G(n,n-k)$.
\end{lem}

We can now formulate the dual counterpart to Theorem
\ref{thm:main2-particular}. Note that since $(L^\circ)^\circ \neq
L$ for a general $k$-Busemann-Petty body, our formulation is a
little different than before.

\begin{thm} \label{thm:main2-dual}
Let $K$ be a centrally symmetric convex body in isotropic
position with $\Vol{K}=\Vol{D_n}$. Then:
\[
L_K \leq C \inf \set{\left . \frac{\L_{2k}^2}{\MR_k(T(L))}  \right
|
\begin{array}{c} L \subset K^\circ \; , \; L \in \BP_k^n \; , \\ T \in SL(n) \; , \; k =
1,\ldots,\lfloor n/3 \rfloor \end{array} }.
\]
\end{thm}
\begin{proof}
First, let us assume $\Vol{K}=1$, and correct for this later. Fix
$k=1,\ldots,\lfloor n/3 \rfloor$ and let $L \in \BP_k^n$ be
contained in $K^\circ$. As in the proof of Theorem
\ref{thm:main1-dual}, we note that $T(L) \in \BP_k^n$ for any $T
\in SL(n)$. Applying Theorem \ref{thm:main2-basic}, the left hand
side of (\ref{eq:main2-basic}) gives:
\begin{equation} \label{eq:main2-polar}
L_K  / C_1 \geq \brac{ \frac{ \widetilde{V}_k(T(L),D) }
{\widetilde{V}_k(T(L),K)}}^{1/k} \geq \brac{ \frac{
\widetilde{V}_k(T(L),D) } {\widetilde{V}_k(T(K^\circ),K)}}^{1/k},
\end{equation}
for $D$ the Euclidean ball of volume 1. Evaluating the denominator
on the right using the trivial $\norm{x}_{T(K^\circ)} \norm{x}_K
\geq \abs{\scalar{T^{-1}(x),x}} = \abs{T^{-1/2}(x)}^2$ for any
positive definite $T\in SL(n)$, we have that:
\begin{eqnarray}
\nonumber & & \brac{ \widetilde{V}_k(T(K^\circ),K) }^{1/k} =
\brac{\frac{1}{n} \int_{S^{n-1}} \norm{x}^{-k}_{T(K^\circ)}
\norm{x}_K^{-(n-k)} dx}^{1/k} \\
\nonumber & \leq & \brac{\frac{1}{n} \int_{S^{n-1}}
\norm{x}_{T^{1/2}(D_n)}^{-2k} \norm{x}_K^{-(n-2k)}dx}^{1/k} =
V_{2k}(T^{1/2}(D_n),K)^{1/k}.
\end{eqnarray}
Using property (\ref{eq:dmv-property}) of dual mixed-volumes, the
latter expression is equal to $V_{2k}(D_n,T^{-1/2}(K))^{1/k}$.
Denoting $G = G(n,n-2k)$, and using polar coordinates and double
integration, we have:
\begin{eqnarray}
\nonumber & & \!\!\!\!\!\!\!\!\!\! V_{2k}(D_n,T^{-1/2}(K))^{1/k} =
\brac{ \Vol{D_n} \int_{G} \int_{S^{n-1}\cap E} \norm{\theta}_{T^{-1/2}(K)}^{-(n-2k)} d\sigma_E(\theta) d\nu(E) }^{1/k} \\
\nonumber & = & \brac{ \frac{\Vol{D_n}}{\Vol{D_{n-2k}}} \int_{G}
\Vol{T^{-1/2}(K) \cap E} d\nu(E) }^{1/k}  \leq
\frac{C}{n-2k}\brac{\frac{\L_{2k}}{L_K}}^2,
\end{eqnarray}
where we have used Lemma \ref{lem:average-sections} in the last
inequality and (\ref{eq:ball-volume-formula}).
Together with (\ref{eq:main2-polar}), cancelling out one $L_K$
term, and using $n - 2k \geq n/3$, this gives:
\begin{equation} \label{eq:main2-dual-L_K-bound}
L_K \leq C' n^{-1} \frac{\L^2_{2k}}{\widetilde{V}_k(T(L),D)^{1/k}}
\simeq n^{-1/2} \frac{\L^2_{2k}}{\MR_k(T(L))},
\end{equation}
for any $T \in SL(n)$ (since it can be factorized into a
composition of a rotation and a positive-definite transformation,
and $\MR_k$ is invariant to rotations). Now correcting for our
initial assumption on $\Vol{K}$ and going back to
$\Vol{K}=\Vol{D_n}$, we have the desired:
\[
L_K \leq C \frac{\L^2_{2k}}{\MR_k(T(L))}.
\]
\end{proof}

As in the previous section, it would be nice to know that for $L
\in \BP_k^n$, there exists a position in which we can bound
$\MR_k(T(L))$ from below by $\sVolRatio{L}{D_n}$ times some
function of $k$. Unfortunately, we cannot provide an analogue of
Lemma \ref{lem:M_p-bound} for general $k$-Busemann-Petty bodies,
but for \emph{convex} members we have the following lemma, which
is stated again in Corollary \ref{cor:bounds-on-mean-radius}:
\begin{lem} \label{lem:MR_k-bound}
Let $K$ be an isotropic convex body with $\Vol{K} = \Vol{D_n}$,
and assume that $K\in \BP_k^n$ for some $k=1,\ldots,n-1$. Then:
\[
\MR_k(K) \geq C / \Lk.
\]
\end{lem}
\begin{proof}
This is a trivial consequence of Theorem
\ref{thm:main2-particular} applied with $L=K$, and using the well
known fact (e.g. \cite{Milman-Pajor-LK}) that $L_K$ is always
bounded from below by a universal constant.
\end{proof}

We will therefore require that the body $L$ from Theorem
\ref{thm:main2-dual} be convex, and denote by $C\BP_k^n$ the class
of convex $k$-Busemann-Petty bodies in $\Real^n$. Applying Lemma
\ref{lem:MR_k-bound} to the body $L$, using the reverse
Blaschke-Santalo inequality (\ref{eq:reverse-Santalo}) and
homogeneity, we immediately have:

\begin{cor}
For any centrally symmetric convex body $K$:
\[
L_K \leq C \inf \set{\left . \L_{2k}^2 \Lk \VolRatio{L^\circ}{K}
\right | \begin{array}{c} K \subset L^\circ \; , \; L \in
C\BP_k^n \; , \\ k = 1,\ldots,\lfloor n/3 \rfloor \end{array} }.
\]
\end{cor}

We will see some applications of Theorem
\ref{thm:main2-particular} in the next section.


\section{Applications} \label{sec:applications}

As applications, we state a couple of immediate consequences of
Corollaries \ref{cor:main1-particular} and \ref{cor:main1-dual}
about the isotropic constant of polytopes with few facets or
vertices. Next, we give several corollaries of Theorem
\ref{thm:main2-particular}, and show how they may be used to
bound the isotropic constant of new classes of bodies.

It is well known that any centrally symmetric polytope with $2m$
facets is a section of an $m$-dimensional cube, and by duality,
any centrally symmetric polytope with $2m$ vertices is a
projection of an $m$-dimensional unit ball of $l_1$. It is also
well known that $l_\infty^m$ isomorphically embeds in $L_p$ for
$p=\log(1+m)$, and by duality, $l_1^m$ is isomorphic to a quotient
of $L_q$, for $q=p^*$ the conjugate exponent to $p$. With the same
notations, it follows that a polytope with $2m$ facets is
isomorphic to a section of $L_{p}$ and that a polytope with $2m$
vertices is isomorphic to a quotient of $L_q$. The following is
therefore an immediate consequence of Corollary
\ref{cor:main1-particular} or Junge's Theorem:

\begin{cor}
Let $K$ be a convex centrally symmetric polytope with $2m$
facets. Then $L_K \leq C \sqrt{\log(1+m)}$.
\end{cor}

\noindent Since any convex body may be isomorphically approximated
by a polytope with $C^n$ facets (or vertices), we retrieve the
well known naive bound $L_K \leq C \sqrt{n}$. In this respect, the
factor of $\sqrt{p}$ in Corollary \ref{cor:main1-particular} for
sections of $L_p$ seems more natural than the factor of $p$ for
quotients of $L_q$, appearing in Corollary \ref{cor:main1-dual} or
Junge's Theorem. Reproducing the above argument, an immediate
consequence of Corollary \ref{cor:main1-dual} or Junge's Theorem
is:

\begin{cor} \label{cor:Junge}
Let $K$ be a convex centrally symmetric polytope with $2m$
vertices. Then $L_K \leq C \log(1+m)$.
\end{cor}

\noindent As mentioned in the Introduction, Corollary
\ref{cor:Junge} implies that Gluskin's probabilistic construction
in \cite{Gluskin-Diameter} of two convex bodies $K_1$ and $K_2$
with Banach-Mazur distance of order $n$, satisfies
$L_{K_1},L_{K_2} \leq C \log(1+n)$. This is simply because the
bodies $K_1$ and $K_2$ are constructed as random polytopes with
(at most) $4n$ vertices.

\medskip

Another easy corollary, which was already partially stated in
Lemmas \ref{lem:M_p-bound} and \ref{lem:MR_k-bound}, may be
deduced from Theorems \ref{thm:L_K-leq-type-2},
\ref{thm:main1-particular} and \ref{thm:main2-particular}, if we
use the well known fact that $L_K$ is always bounded from below.
Together with Jensen's inequality (as in (\ref{eq:Jensen})), this
reads as follows:

\begin{cor} \label{cor:bounds-on-mean-radius}
Let $K$ be convex centrally symmetric isotropic body with
$\Vol{K}=\Vol{D_n}$. Then:
\begin{enumerate}
\item
$1 \leq M_2(K) \leq C T_2(X_K)$.
\item
If $K \in SL_p^n$ ($p\geq 0$), then $1 \leq M_p(K) \leq C
\sqrt{p_0}$, where $p_0 = \max(1,\min(p,n))$.
\item
If $K \in \BP_k^n$ ($k=1,\ldots,n-1$), then $C / \Lk \leq \MR_k(K)
\leq 1$.
\end{enumerate}

\end{cor}

\medskip

Next, we proceed to deduce several consequences of Theorem
\ref{thm:main2-particular}. It is known that $\BP_k^n$ does not
contain all convex bodies for $k<n-3$, and that $\BP_{n-1}^n$
already contains all star-bodies
(\cite{Bourgain-Zhang},\cite{Koldobsky-I-equal-BP}). So definitely
not all convex bodies are isometric to members of $\BP_k^n$ for
$k<n-3$. Nevertheless, the following assumption might be true:

\medskip
\noindent\textbf{Outer Volume Ratio Assumption for $\BP_k^n$.}
\emph{There exist two universal constants $C,\epsilon>0$, such
that for any $n$ and any convex body $K$ in $\Real^n$ there exists
a star-body $L \in \BP_k^n$ for $k = n^{1-\epsilon}$, such that $K
\subset L$ and $(\Vol{L}/\Vol{K})^{1/n} \leq C$.}

\smallskip
Under this assumption, Theorem \ref{thm:main2-particular} would
immediately imply that $\L_n \leq C \L_{n^{1-\epsilon}}$.
Denoting $\delta = -1 / \log(1-\epsilon)$, and iterating this
inequality $\delta \log \log n$ times, we would have:

\begin{cor}
Under the Outer Volume Ratio Assumption for $\BP_k^n$, we have
$\L_n \leq C_1 (\log(1+n))^{C_2 \delta}$ for $\delta>0$ as above.
\end{cor}

\smallskip
In addition, the advantage of working with $\BP_k^n$ when trying
to find or build a body $L \in \BP_k^n$ containing $K$, is that
we need not worry about the convexity of $L$ like in the case of
$SL_p^n$. The convexity of $K$ has already been used in Theorem
\ref{thm:main2-basic} (in (\ref{eq:sections-of-isotropic})), so we
may now consider $\rho_K$ as a function on $S^{n-1}$ which we
want to tightly bound from above using functions $\rho_L$ from
the given family $\BP_k^n$. This is an especially attractive
approach, as $\BP_k^n$ has the following nice characterization,
first proved by Goodey and Weil in \cite{Goodey-Weil} for
intersection-bodies (the case $k=1$), and extended to general $k$
by Grinberg and Zhang in \cite{Grinberg-Zhang}:

\smallskip
\noindent\textbf{Theorem (Grinberg and Zhang). }\emph{A star-body
$K$ is a $k$-Busemann-Petty body iff it is the limit
of $\set{K_i}$ in the radial metric $d_r$, where each $K_i$ is a
finite $k$-radial sums of ellipsoids $\set{\E^i_j}$:
\[
\rho^k_{K_i} = \rho^k_{\E^i_1} + \ldots + \rho^k_{\E^i_{m_i}},
\]
or equivalently, if there exists a Borel measure $\mu$ on $SL(n)$
such that:
\[
\rho^k_{K_i} = \int_{SL(n)} \rho^k_{T(D_n)} d\mu(T).
\]
}

In fact, even the "easiest" case $k=1$ in Theorem
\ref{thm:main2-particular} seems potentially useful, as we shall
demonstrate below. Note that since the intersection-body $L$ need
not be convex (and therefore Corollary
\ref{cor:bounds-on-mean-radius} does not apply to it), the
mean-radius $\MR(L)$ might be significantly smaller than the
volume-radius $(\Vol{L}/\Vol{D_n})^{1/n}$. As demonstrated by
Theorem \ref{thm:main2-particular}, a smart way to bound $\rho_K$
from above by $\rho_L$ which is the sum of radial functions of
ellipsoids, such that we have control over $L$'s mean-radius,
might provide a new bound on the isotropic constant. We give two
examples of how such an approach might work. Unfortunately, we
need to use some additional assumptions, which, although we
believe to be true, we have not been able to prove. First, we
need a new definition for a class of bodies.

\definition{
Let $K$ denote a star-body. We will work with the radial metric
topology on the space of star-bodies. Introduce the closed set of
volume preserving linear images of $K$,
\[B(K) = \set{T(K) \; | \; T \in SL(n)}.\]
The \emph{Radial Sums of $K$}, denoted by $RS(K)$, is the
closure in the radial metric of the
family of all star-bodies $L$, such that there exists a
non-negative Borel measure $\mu$ on $B(K)$, for which:
\[
\rho_L = \int_{B(K)} \rho_{K'} d\mu(K').
\]

\noindent Similarly, if $P$ is a closed set of star-bodies, then
the \emph{Radial Sums of $P$}, denoted $RS(P)$, is the
closure in the radial metric of the
family of all star-bodies $L$, such that there exists a
non-negative Borel measure $\mu$ on $B(P) = \bigcup_{K \in P}
B(K)$, for which:
\[
\rho_L = \int_{B(P)} \rho_{K'}(\theta) d\mu(K').
\]
}

So for example $RS(D_n)$ is exactly the class of
intersection-bodies, since $B(D_n)$ is the set of all ellipsoids
of volume $\Vol{D_n}$, and by the aforementioned result of Goodey
and Weil, the radial sums of this set are exactly the class of
intersection-bodies. Another easy observation is that $RS(P)$ is
closed under full-rank linear transformations, since for any
linear $T$:
\[
\rho_K = \rho_{K_1} + \rho_{K_2} \; \Rightarrow \; \rho_{T(K)} =
\rho_{T(K_1)} + \rho_{T(K_2)}.
\]
As a consequence, $RS(D_n) \subset RS(K)$ for any star-body $K$.
To see this, first notice that $D_n \in RS(K)$, by choosing the
Borel measure $\mu$ on $B(K)$ to be:
\[
\mu(A) = \eta(\set{T \in O(n) \; | \; T(K) \in A})
\]
for every Borel set $A \subset B(K)$, where $\eta$ is the
appropriately normalized Haar measure on $O(n)$, the group of
orthogonal rotations in $\Real^n$. Since $RS(K)$ is closed under
$SL(n)$, radial summation, and limit in the radial-metric, it
follows that $RS(D_n) \subset RS(K)$. Therefore, for any non
intersection-body $K$, $RS(K)$ properly contains the class of
intersection bodies.

There are many interesting questions that may be asked about
Radial Sums of star-bodies, such as whether it is possible to
characterize a minimal set $P$ for which $RS(P)$ already contains
all convex bodies, or, probably easier, all star-bodies. In
particular it is not even clear to us whether $P$ may be chosen
as a singleton in either case. Our focus will be on the following
two assumptions, which we believe to be true. The first is about
the $n$-dimensional cube $Q_n$ (of volume 1):

\medskip
\noindent\textbf{Outer Mean-Radius Assumption for the Cube $Q_n$.}
\emph{For any $K\in B(Q_n)$, there exists an ellipsoid $\E$
containing $K$ such that $\MR(\E) / \MR(K) \leq C \log(1+n)$, for
some universal constant $C>0$.}
\medskip

The second assumption is about $UC(n)$, the class of volume 1
convex bodies in $\Real^n$ which are all unconditional with
respect to the same fixed Euclidean structure. We shall say that
a body is a cross-polytope if it is a linear-image of the unit
ball of $l^n_1$.

\medskip
\noindent\textbf{Outer Mean-Radius Assumption for UC(n).}
\emph{For any $K \in B(UC(n))$, there exists a cross-polytope $L$
containing $K$ such that $\MR(L) / \MR(K) \leq C \log(1+n)$, for
some universal constant $C>0$. }
\medskip

We will shortly give motivation for why these assumptions might
be correct, but first, let us show an easy consequence of Theorem
\ref{thm:main2-particular} under each assumption.

\begin{cor} \hfill
\begin{enumerate}
\item
Under the Outer Mean-Radius Assumption for $Q_n$, for any convex
body $K \in RS(Q_n)$, we have $L_K \leq C \log(1+n)$.
\item
Under the Outer Mean-Radius Assumption for $UC(n)$, for any convex
body $K \in RS(UC(n))$, we have $L_K \leq C \log(1+n)$.
\end{enumerate}
\end{cor}
\noindent As mentioned before, the families of convex bodies in
$RS(Q_n)$ and $RS(UC(n))$ are potentially new classes of convex
bodies, which might contain a big piece of the convex bodies
compactum. Therefore, this new approach to bounding the isotropic
constant might be applicable for a large family of convex bodies.
\begin{proof}
Let $K$ be an isotropic convex body of volume $\Vol{D_n}$ in
$RS(P)$, where $P$ is either $\set{Q_n}$ or $UC(n)$. By
approximation, we may assume that $\rho_K = \sum_i \mu_i
\rho_{K_i}$, where $K_i \in B(P)$ and $\mu_i \geq 0$.

Notice that both the unit-ball of $l_1^n$ and the Euclidean ball
are intersection bodies, and this is preserved under volume
preserving linear transformations. Therefore, by the Outer
Mean-Radius Assumption for $P$, there exist intersection-bodies
$L_i$ such that $K_i \subset L_i$ and $\MR(L_i) / \MR(K_i) \leq C
\log(1+n)$. Now define $L$ to be the star-body for which $\rho_L
= \sum_i \mu_i \rho_{L_i}$. It is obvious that $L$ contains $K$,
and that $L$ is an intersection-body (since these are closed
under non-negative radial summation, as follows from their
definition). In addition, since the mean-radius $\MR$ is additive
under radial summation, it is clear that $\MR(L) / \MR(K) \leq C
\log(1+n)$. But using Jensen's inequality (as in
(\ref{eq:Jensen})), we have $\MR(K) \leq \MR_n(K) = 1$, and
therefore $\MR(L) \leq C \log(1+n)$. Using Theorem
\ref{thm:main2-particular}, the proof is complete.
\end{proof}

We conclude by giving motivation for why the above two
assumptions might be correct, and explain the difficulty in
proving them. The next proposition demonstrates that the
assumptions indeed hold when the bodies in question are in
isotropic position, in which case the bounding bodies may be
chosen to be in isotropic position as well.

\begin{prop} \label{prop:assumptions-work-for-isotropic}
\hfill
\begin{enumerate}
\item
Let $D$ be the circumscribing Euclidean ball of $Q_n$. Then:
\[\MR(D) / \MR(Q_n) \leq C \log(1+n).\]
\item
Let $K$ be an unconditional convex body in isotropic position,
and let $L$ be its circumscribing unit ball of $l_1^n$. Then:
\[\MR(L) / \MR(K) \leq C \log(1+n).\]
\end{enumerate}
\end{prop}

\begin{proof}
\hfill
\begin{enumerate}
\item
This is a standard calculation relating to the concentration of
the norm $\norm{\cdot}_{Q_n}$ on the sphere, which may be done
using the standard concentration techniques from
\cite{Milman-Schechtman-Book}.
We prefer to quote a general result by Klartag and Vershynin from
\cite[Proposition 1.2]{Klartag-Vershynin}, which states that for
any convex body $K$, if $0<l<C k(K)$, where $k(K) =
n(M(K)/b(K))^2$, then $\MR_l(K) \simeq 1/M(K)$. Since for the
volume 1 cube $Q_n$ it is well known (e.g.
\cite{Milman-Schechtman-Book}) that $M(Q_n) \simeq
\sqrt{\log(1+n)} / \sqrt{n}$, $b(Q_n) = 2$, and therefore $k(Q_n)
\simeq \sqrt{\log(1+n)}$, it follows that for $n$ large enough we
may use the above result for $l=1 < C k(Q_n)$, to conclude that
(for all $n$) $\MR(Q_n) \simeq \sqrt{n} / \sqrt{\log(1+n)}$. Since
$\MR(D) = \sqrt{n} / 2$, the claim follows.

\item
Let $P_n$ be the unit ball of $l_1^n$ of volume 1. It is well
known (e.g. \cite{Bobkov-Nazarov}) that there exist $C_1,C_2
> 0$, such that for any isotropic convex body $K$ of volume 1,
which is unconditional with respect to the given Euclidean
structure, the following inclusions hold:
\[
C_1 Q_n \subset K \subset C_2 P_n. \] Therefore $\MR(L) / \MR(K)
\leq \MR(C_2 P_n) / \MR( C_1 Q_n)$. We have already seen that
$\MR(Q_n) \simeq \sqrt{n} / \sqrt{\log(1+n)}$. We may estimate
$\MR(P_n)$ in the same manner, or alternatively, use Corollary
\ref{cor:bounds-on-mean-radius} to deduce that $\MR(P_n) \simeq
\sqrt{n}$. Therefore $\MR(L) / \MR(K) \leq C \log(1+n)$.
\end{enumerate}
\end{proof}

Unfortunately, the techniques described above fail when used upon
$T(K)$, where $K$ is in isotropic position but $T$ is an almost
degenerate mapping. In particular, it is a bad idea to try to
bound $T(K)$ using $T(L)$, where $L$ is the optimal bounding body
for $K$. Indeed, let us try to evaluate $\MR(T(D))/\MR(T(Q_n))$,
where as in Proposition
\ref{prop:assumptions-work-for-isotropic}, $D$ is the
circumscribing Euclidean ball of $Q_n$. Using
(\ref{eq:dmv-property}), we have:
\[
\frac{\MR(T(D))}{\MR(T(Q_n))} =
\frac{\widetilde{V}_1(T(D),D_n)}{\widetilde{V}_1(T(Q_n),D_n)} =
\frac{\widetilde{V}_1(D,T^{-1}(D_n))}{\widetilde{V}_1(Q_n,T^{-1}(D_n))}.
\]
Denoting $\E = T^{-1}(D_n)$, we see that:
\[
\frac{\MR(T(D))}{\MR(T(Q_n))} = \frac{\int_{S^{n-1}}
\rho_D(\theta) \rho_\E(\theta)^{n-1}
d\sigma(\theta)}{\int_{S^{n-1}} \rho_{Q_n}(\theta)
\rho_\E(\theta)^{n-1} d\sigma(\theta)},
\]
and this is clearly invariant under homothety of $\E$. Now let us
define $\E(\xi,a,b)$ for $\xi \in S^{n-1}$ and $a,b > 0$ as the
ellipsoid whose corresponding norm is defined as:
\[
\norm{x}^2_{\E(\xi,a,b)} = \frac{\scalar{x,\xi}^2}{a^2} +
\frac{\abs{x}^2 - \scalar{x,\xi}^2}{b^2}.
\]
It was shown in \cite{Goodey-Weil} that by appropriately choosing
$a=a(\epsilon)$ very large and $b=b(\epsilon)$ very small, and
setting $\E(\xi,\epsilon) = \E(\xi,a(\epsilon),b(\epsilon))$, the
family $\rho^{n-1}_{\E(\xi,\epsilon)}$ is an approximation of
unity on $S^{n-1}$ at $\xi$ (as $\epsilon>0$ tends to 0). This
means that for every $f \in C(S^{n-1})$:
\[
\int_{S^{n-1}} f(\theta) \rho^{n-1}_{\E(\xi,\epsilon)}(\theta)
d\sigma(\theta) \longrightarrow f(\xi) \; \textrm{ as } \epsilon
\rightarrow 0.
\]
Hence, we see that by choosing $T = T(\xi)$ to be very degenerate,
we may arbitrarily approximate:
\[
\frac{\MR(T(D))}{\MR(T(Q_n))} \simeq
\frac{\rho_D(\xi)}{\rho_{Q_n}(\xi)},
\]
and the latter ratio may be chosen to be any number between $1$
and $\sqrt{n}$ by an appropriate choice of $\xi \in S^{n-1}$.
This example demonstrates the difficulty in proving the Outer
Mean-Radius Assumptions.


\bibliographystyle{amsalpha}
\bibliography{../../ConvexBib}

\end{document}